\documentclass[reqno,11pt]{amsart}
\usepackage{amsmath, latexsym, amsfonts, amssymb, amsthm, amscd}
\usepackage{sidecap}
\usepackage{multirow,stmaryrd}
\usepackage{color}
\usepackage[colorlinks=true, pdfstartview=FitV, linkcolor=blue, citecolor=blue, urlcolor=blue,pagebackref=false]{hyperref}
\usepackage{graphicx}
\numberwithin{equation}{section}

\newtheorem{theorem}{Theorem}[section]
\newtheorem{lemma}[theorem]{Lemma}

\newcommand{\ind}{\mathbf{1}}

\renewcommand{\tilde}{\widetilde}

\newcommand{\bP}{{\ensuremath{\mathbf P}} }
\newcommand{\bE}{{\ensuremath{\mathbf E}} }

\DeclareMathSymbol{\leqslant}{\mathalpha}{AMSa}{"36} % nicer `smaller or equal'
\DeclareMathSymbol{\geqslant}{\mathalpha}{AMSa}{"3E} % nicer `larger or equal'
\DeclareMathSymbol{\eset}{\mathalpha}{AMSb}{"3F}     % nicer `emptyset'
\renewcommand{\leq}{\;\leqslant\;}                   % redef. of < or =
\renewcommand{\geq}{\;\geqslant\;}                   % redef. of > or =
\newcommand{\sumtwo}[2]{\sum_{\substack{#1 \\ #2}}} % sum with 2 lines
\newcommand{\limtwo}[2]{\lim_{\substack{#1 \\ #2}}}     % \lim with 2 lines
\newcommand{\bbN}{{\ensuremath{\mathbb N}} }

\newcommand{\bbR}{{\ensuremath{\mathbb R}} }

\newcommand{\bbZ}{{\ensuremath{\mathbb Z}} }

\newcommand{\ga}{\alpha}

\newcommand{\gd}{\delta}
\newcommand{\gep}{\varepsilon}       % \ge already exists...

\newcommand{\gl}{\lambda}

\newcommand{\gs}{\sigma}

\makeatletter
\def\captionfont@{\footnotesize}
\def\captionheadfont@{\scshape}

\long\def\@makecaption#1#2{%
  \vspace{2mm}
  \setbox\@tempboxa\vbox{\color@setgroup
    \advance\hsize-6pc\noindent
    \captionfont@\captionheadfont@#1\@xp\@ifnotempty\@xp
        {\@cdr#2\@nil}{.\captionfont@\upshape\enspace#2}%
    \unskip\kern-6pc\par
    \global\setbox\@ne\lastbox\color@endgroup}%
  \ifhbox\@ne % the normal case
    \setbox\@ne\hbox{\unhbox\@ne\unskip\unskip\unpenalty\unkern}%
  \fi
  \ifdim\wd\@tempboxa=\z@ % this means caption will fit on one line
    \setbox\@ne\hbox to\columnwidth{\hss\kern-6pc\box\@ne\hss}%
  \else % tempboxa contained more than one line
    \setbox\@ne\vbox{\unvbox\@tempboxa\parskip\z@skip
        \noindent\unhbox\@ne\advance\hsize-6pc\par}%
\fi
  \ifnum\@tempcnta<64 % if the float IS a figure...
    \addvspace\abovecaptionskip
    \moveright 3pc\box\@ne
  \else % if the float IS NOT a figure...
    \moveright 3pc\box\@ne
    \nobreak
    \vskip\belowcaptionskip
  \fi
\relax
}
\makeatother
\def\writefig#1 #2 #3 {\rlap{\kern #1 truecm
\raise #2 truecm \hbox{#3}}}

\newcommand{\tn}{\textsc{n}}
\newcommand{\tf}{\textsc{f}}

\newcommand{\cM}{\mathcal{M}}

\newcommand{\htx}{\hat\tau^{(1)}}
\newcommand{\hty}{\hat\tau^{(2)}}
\newcommand{\hP}{\bP}

\title[DNA melting structures in the $\mathrm{gPS}$ model]{
DNA melting 
structures \\
in the generalized Poland-Scheraga model
}

\author{Quentin Berger}
\address{
  UPMC, Laboratoire de Probabilit{\'e}s et Mod\`eles Al\'eatoires, UMR 7599,
            F- 75205 Paris,France
}

\author{Giambattista Giacomin}
\address{
  Universit\'e Paris Diderot, Sorbonne Paris Cit\'e,  Laboratoire de Probabilit{\'e}s et Mod\`eles Al\'eatoires, UMR 7599,
            F- 75205 Paris, France
}

\author{Maha Khatib}
\address{
Universit\'e Libanaise, Laboratoire de Math\'ematiques-EDST, Beyrouth, Liban 
            }
\begin{document}

\begin{abstract}
The Poland-Scheraga model for DNA denaturation, besides playing a central role in applications, has been   widely studied  in the physical and mathematical literature over the past decades.
%In the 2000s,
More recently a natural generalization has been introduced in the biophysics literature (see in particular \cite{cf:GO}) to overcome the limits of the original model, namely to allow an {\sl excess of bases} -- i.e. a different length of the two single stranded DNA chains --
 and to allow slippages in the chain pairing. The increased complexity of the model is reflected in the appearance of configurational transitions when the DNA is in  double stranded form.  In \cite{cf:GK} the generalized Poland-Scheraga model has been analyzed thanks to a representation in terms of a bivariate renewal process. In this work we exploit  this representation farther 
 and fully characterize the path properties of the system, making therefore explicit the geometric structures -- and the configurational transitions -- that are observed when the polymer is in the double stranded form. 
What we prove is that when the excess of bases is not   absorbed in a homogeneous fashion  along the double stranded chain -- a case treated in 
\cite{cf:GK} -- then it either condensates in a single  {\sl macroscopic} loop or it accumulates into an unbound single strand free end.
 %More precisely, we show that in the so-called non-Cram\'er regime, only two configurations can occur: either there is a single macroscopic loop and the unbound free ends have length of order one, or there is a macroscopic unbound free end and no macroscopical loop -- in any case, all non-macroscopical loops are of smaller (polynomial) order. This is in contrast with the Cram\'er regime, where the unbound free ends have length of order one, and all the loops are (logarithmically) small.
   \\[8pt]
  \\[8pt]
\end{abstract}

\maketitle

%\tableofcontents

\section{Introduction}
\label{sec:intro}

The
Poland-Scheraga (PS) model \cite{cf:PSbook,cf:BD,cf:EON} played and still plays a central role in the analysis of DNA denaturation (or melting): double stranded DNA  \emph{melts} into two single stranded DNA polymer chains at high temperature. The success of the model is  partly due to the fact that it is exactly solvable  when the heterogeneous character of the DNA 
%-- 
%C-G bonds are stringer than A-T bonds 
%-- 
is neglected. Moreover, solvability has an interest on its own, from a more theoretical standpoint: phase transition and critical phenomena in the PS model are completely understood \cite{cf:Fisher,cf:GB}. However, the PS model is an oversimplification in many respects: it deals with two strands of equal length and it does not allow slippages of the two chains. These  simplifications make the model one dimensional, and solvability becomes less surprising. What is instead surprising is that a natural generalization \cite{cf:GO0, cf:GO,cf:NG} -- called generalized Poland-Scheraga (gPS) model -- fully overcomes these limitations,  retaining the solvable character in spite of the substantially richer variety of structures that it displays.
In  \cite{cf:GK} a mathematical approach to the gPS model is developed and it is pointed out that
it can be represented 
 in terms of a two dimensional renewal process, much like the PS model can be represented in terms of a one dimensional renewal. The solvable character of both models is then directly related to their renewal structure.
The growth in complexity from PS to gPS models is nevertheless considerable: 
the key feature of PS and gPS  is the presence of a localization transition, corresponding to the passage from separated to bound strands, and  for the gPS there are three, not only one, types of localized trajectories (or configurations). 
This has been first pointed out, at least in part, in \cite{cf:NG}, where one can find theoretical arguments (based also on a Bose-Einstein condensation analogy) and numerical evidence  that \emph{``suggest that a temperature-driven conformational transition occurs before the melting transition"} \cite[p.3]{cf:NG}. 
 %favor of the following phenomenology:

In this work we fully characterize the possible localized configurations.  
%(there can be more than one \cite{cf:GK}). 
The transitions between different types of configurations  have been already studied at the level of the free energy in \cite{cf:GK} where  these phenomena have been  mathematically identified and  interpreted in a Large Deviations framework  in terms of \emph{Cram\'er} and \emph{non-Cram\'er} strategies. This will be explained in detail below. Here we content ourselves with pointing out that a full analysis of the Cram\'er regimes is given in \cite{cf:GK}. However, the non-Cram\'er regime, where the condensation phenomena happen,  requires a a substantially finer analysis -- moderate deviations and local limit estimates -- at the level of the bivariate renewals.
These estimates, to which much attention has been devoted in the literature in the one dimensional set-up (see \cite{cf:AL,cf:DDS08} and references therein),  are lacking to the best of our knowledge for higher dimensional renewals and they are not straightforward generalizations. They represent 
the technical core of this paper.

\subsection{The Model and some basic results}
\label{sec:model}

\begin{figure}[htbp]
\centering
\includegraphics[width=11 cm]{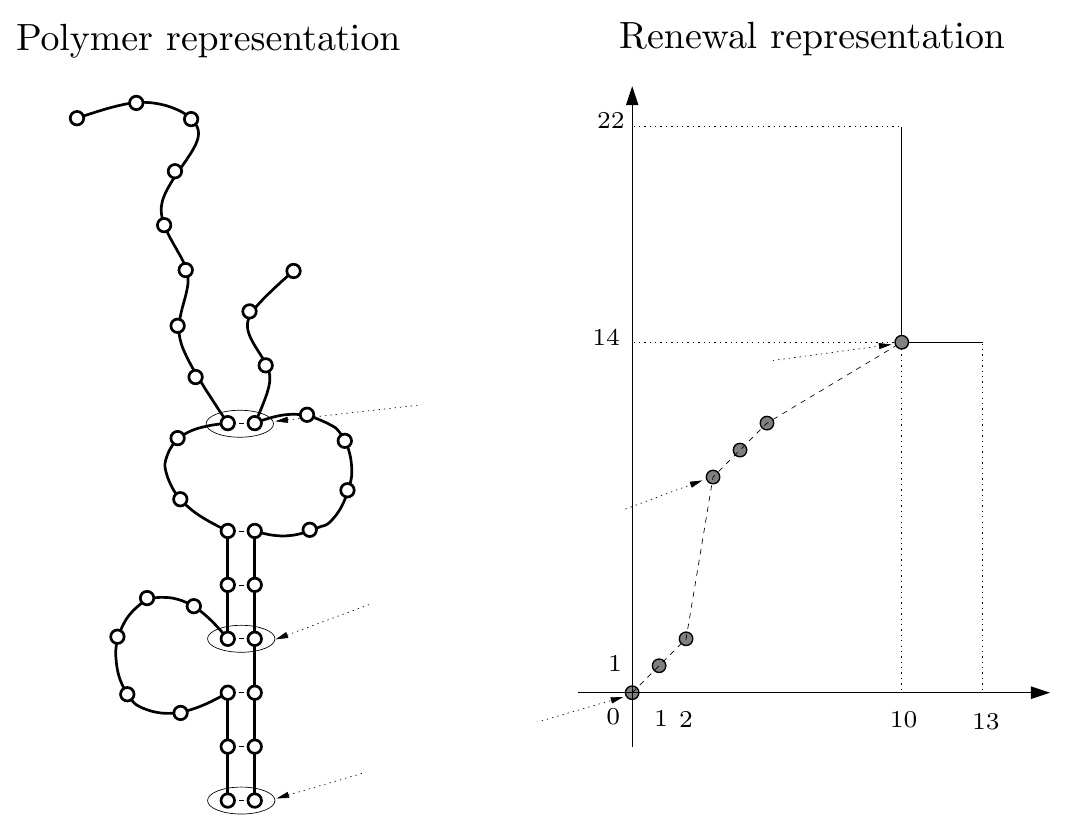}
\vskip-.2cm
\caption{\label{fig:explain} 
A configuration of the gPS model, with one strand containing 23 bases and the other 14,
is represented in two fashions: the \emph{natural} (or polymer) one and the renewal one.  In particular  we see that $(1,1)$ renewal increments (or inter-arrivals) correspond to bound base pairs and all other increments $(i,j)$ correspond to unbound regions in the bulk, that we call \emph{loops} (of total length $i+j$, with length $i$ in the first strand and $j$ in the second strand). The term \emph{unbound} is rather reserved to the terminal portion of the polymer: we refer the free ends as unbound strands. 
Throughout this work, a polymer trajectory is always given in the renewal representation: it is therefore just a point process in the plane. 
}
\end{figure}

We introduce the model in detail only from the renewal representation. The link with the original representation  of the model is summed up
in Fig.~\ref{fig:explain} and its caption, and we refer to \cite{cf:GK} for more details.

We consider a persistent bivariate renewal process $\tau = { \lbrace (\tau_n^{(1)},\tau_n^{(2)}) \rbrace }_{n \ge 0}$, that is a sequence of random variables such that $\tau_0 =(0,0)$, ${ \lbrace \tau_n - \tau_{n-1} \rbrace}_{n=0,2, \ldots}$ is IID and such that the inter-arrival law -- {i.e.}\ the law of $\tau_1$ --, takes values in $\bbN^2 := \lbrace 1, 2, \cdots \rbrace ^2$.

We set $\bP ( \tau_1 = (n,m))= K(n+m)$ with 
\begin{equation}
K(n) := \frac{L(n)}{n^{2 + \ga}} \,,
\end{equation}
for some $\ga > 0$ and some slowly varying function $L(\cdot)$. Moreover $\sum_{n,m} K(n+m) = 1$ since we assumed the process to be persistent.

We consider two versions of the model: constrained and free.
The  partition function of the constrained model, or \emph{constrained partition function}, can be written as
\begin{equation}
Z^c_{N,M,h} := \sum_{n=1}^{N \wedge M} \sumtwo{\underline{l} \in \bbN^n}{\vert \underline{l} \vert = N} \sumtwo{\underline{t} \in \bbN^n}{\vert \underline{t} \vert = M} \prod_{i=1}^N \exp(h) K(l_i+t_i) \, ,
\end{equation}
where $h\in\bbR$ is the binding energy, or pinning parameter.

The partition function of the free model, or \emph{free partition function},  is defined by
\begin{equation}
\label{eq:zf}
Z^f_{N,M,h} := \sum_{i=0}^N \sum_{j=0}^M K_f(i) K_f(j) Z^c_{N-i,M-j,h} \,,
\end{equation}
where $K_f(n) := \overline{L}(n) n^{- \overline{\ga}}$ for some $\overline{\ga} \in \bbR$ and slowly varying function $\overline{L}(\cdot)$. We assume that $K_f(0)=1$ just to prevent this constant from popping up in various formulas: this choice has the side effect of making clear that $K_f(\cdot)$ is not a probability.

In  \cite{cf:GK} it is shown that  for every $h$ and every $\gamma>0$
\begin{equation}
\label{eq:fcf}
\tf_\gamma(h) := \limtwo{N,M \to \infty}{M/N \to \gamma} \frac{1}{N} \log Z^c_{N,M,h} = \limtwo{N,M \to \infty}{M/N \to \gamma} \frac{1}{N} \log Z^f_{N,M,h} < \infty\,,
\end{equation}
which says that the  free energy (density) of free and constrained  models, with binding energy $h$ and strand length
asymptotic ratio equal to $\gamma$, coincide.
A number of basic properties of $h\mapsto \tf_\gamma(h)$ are easily established, notably that it is a convex non decreasing function, equal to zero for $h\le 0$ and positive for $h>0$. This already establishes that $h=0$ is a critical point, in the sense that
$\tf_\gamma(\cdot)$ is not analytic at the origin.

\begin{figure}[htb] 
\centering
\hspace{2cm}
   \begin{minipage}[b]{.35\linewidth}
    \scalebox{0.35}{\includegraphics*{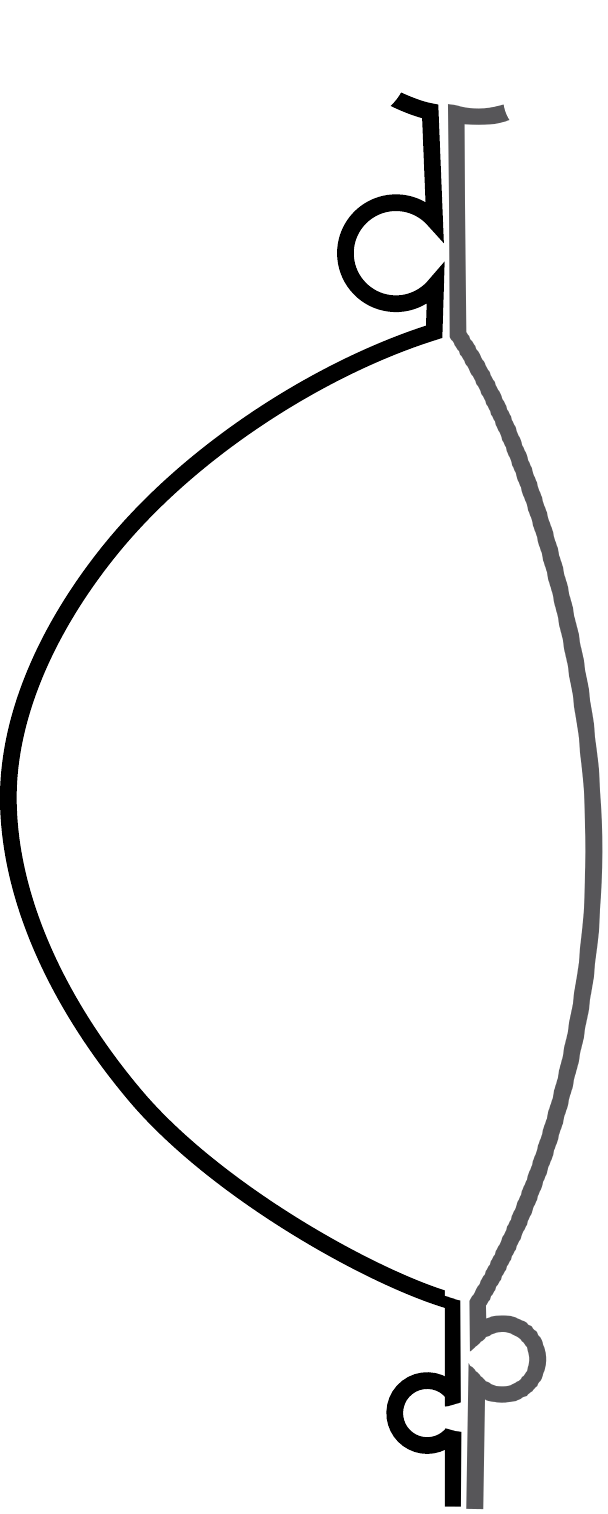}}
   \end{minipage}
    \qquad  
 \begin{minipage}[b]{.35\linewidth}
  \scalebox{0.32}{\includegraphics*{1d2-eps-converted-to.pdf}}
 \end{minipage}\hfill 
\caption{\label{fig:delocalized} A schematic image of the two types of observable trajectories of the free gPS model in the delocalized (denaturated, melted) regime, according to whether the exponent $\overline{\ga} $ is smaller (left picture) or larger (right picture) than $1+ \alpha/2$. In the constrained case only the trajectory on the right is observed, and the small free tails are reduced to zero. This case is treated in \cite{cf:GK}.}
\end{figure}

But \cite{cf:GK} is not limited to results on the free energy:  
associated to $Z^c_{N,M,h}$ and $Z^f_{N,M,h}$ there are two probability measures, that we denote respectively by $\bP_{N, M, h}^c$ and $\bP_{N, M, h}^f$.
They are  point measures, like the renewal processes on which they are built. It is standard to see that
$\partial_h \tf_\gamma (h)$ (which exists except possibly for countably many values of $h$) yields the $N \to \infty$ limit of 
the expected  density of points (under $\bP_{N, M, h}^c$ or $\bP_{N, M, h}^f$). Hence for $h<0$ the density is zero, while 
for $h>0$ the density is positive. This tells us that we are stepping from a regime in which the two strands are essentially fully unbound to a regime in which they are tightly bound. 
In  \cite{cf:GK} results go well beyond this: it is in particular proven that for $h<0$ the number of renewal points is $O(1)$ and these
points are all close to $(0,0)$ or $(N,M)$ (see Fig.~\ref{fig:delocalized}). 
In the polymer representation, this means that the two DNA strands are completely unbound, except for a few contacts between the bases 
just close to the extremities. 
More precisely, it was found in \cite{cf:GK} that in the free case, if $\overline{\ga}  < 1+ \ga/2$ the two strands are free except for $O(1)$ contacts close to the origin, and if $\overline{\ga}  > 1+ \ga/2$ the two free ends are of length $O(1)$ and a large loop appears in the system, see Fig~\ref{fig:delocalized}. 

On the other hand for $h>0$ the situation is radically different. This has has been  analyzed in \cite{cf:GK} but only in the Cram\'er regime. We are now going to discuss this in details.

\subsection{Binding strategies}
\label{sec:strategies}
A way to get a grip on what is going on for $h>0$ is to observe that we can make the elementary manipulation: 
for every non negative $\gl_1$ and $\gl_2$
\begin{equation}
Z^c_{N,M,h} = e^{N \gl_1+ M \gl_2}
\sum_{n=1}^{N \wedge M} \sumtwo{\underline{l} \in \bbN^n}{\vert \underline{l} \vert = N} \sumtwo{\underline{t} \in \bbN^n}{\vert \underline{t} \vert = M} \prod_{i=1}^N \exp\left(h - \gl_1 l_i -\gl_2 t_i\right) K(l_i+t_i) \,. 
\end{equation}
Since $h>0$ we identify a family, in fact a curve in $[0,\infty)^2$, of values of $(\gl_1, \gl_2)$ such that
\begin{equation}
\label{eq:norm1}
\sum_{l, t=1}^\infty
 \exp\left(h - \gl_1 l -\gl_2 t\right) K(l+t) \, =\, 1\, ,
\end{equation}
and \eqref{eq:norm1} clearly defines a probability distribution that is an inter-arrival distribution for a new bivariate renewal process. 
At this point is not too difficult to get convinced that 
$Z^c_{N,M,h}$ is equal to $e^{N \gl_1+ M \gl_2}$ times the probability that this new renewal hits $(N, M)$ (we call this probability {\sl target probability}).
If we are able to choose  $(\gl_1, \gl_2)$ so that the logarithm of the target probability is $o(N)$, then of course $\tf_\gamma (h)=\gl_1+\gamma \gl_2$. This can actually be done: it amounts to solving a variational problem and the uniqueness of the optimal $(\gl_1, \gl_2)$ follows by convexity arguments. However the solution may be qualitatively different
for different values of $h$:

\begin{enumerate}
\item the optimal $(\gl_1, \gl_2)$ belong to $(0, \infty)^2$, so both components of the inter-arrival law of the arising renewal 
have distributions that decay exponentially. We call this \emph{Cram\'er regime} because the tilt of the measure (in both components) is efficient in targeting the point $(N, M)$ to which we are aiming at;
\item either $\gl_1$ or $\gl_2$ is zero, so only one component of the arising inter-arrival law is exponentially tight.
For the sake of conciseness we  call this for now  \emph{non-Cram\'er regime} because the tilt of the measure (in only one of the component) is only partially successful in targeting the point $(N,M)$. 
To be precise there is in reality a boundary region between the two regimes, and the notion of non-Cram\'er regime
will be made more precise just below -- this regime is the main issue of this work -- so we will not dwell further on this right now.
\end{enumerate}

A full treatment of the Cram\'er regime is given in \cite{cf:GK}, and the results can be resumed as follows: all loops are small, in fact the largest is $O( \log N)$, and the unbound strands are of length $O(1)$ -- see the leftmost case in  Fig~\ref{fig:localized}. 
\smallskip

In this work, we focus on the non-Cram\'er regime and the reader who wants to have an anticipation on the results 
can have a look Fig~\ref{fig:localized}.

\subsection{The non-Cram\'er regime}
\label{sec:nonCramer}

In order to make as explicit as possible for which values of $h>0$ the system is in the non-Cram\'er regime, 
let us define $\tn(h)>0$ as the unique solution of
\begin{equation}
\label{eq:Nh}
\sum_{n,m=1}^\infty 
K(n+m) \exp\left(- n \tn(h)\right)\, =\, \exp(-h)\, .
\end{equation}
This computation amounts to solving the variational problem  we were after, in the case in which the problem
is not solvable in $(0, \infty)^2$ and the optimal tilt of the measure involves only one of the two components.  From \eqref{eq:Nh} one can extract a number of properties of $\tn(\cdot)$: it is a real  analytic,  positive,  convex and increasing function \cite{cf:GK}.
We insist on the fact that, in spite of being defined for every $h>0$,  $\tn(h)$ is not always equal to the free energy.
More precisely in \cite{cf:GK} it is shown that
 $\tf_\gamma(h)=
\tn(h)$   if and only if $\gamma \notin (1/\gamma_c(h), \gamma_c(h) )$, where
%$\gamma_c(\cdot): (0,\infty) \longmapsto (1, \infty)$ is a real analytic function and 
\begin{equation}
\label{eq:gammac}
\gamma_c(h) \, :=\, \frac{\sum_{n,m} m K(n+m) \exp(-n \tn(h))}{\sum_{n,m} n K(n+m) \exp(-n \tn(h))}\,.
\end{equation} 
We refer to \cite{cf:GK} for more details on the form of the function $\gamma_c(\cdot)$ and the switching phenomena between the Cram\'er and the non-Cram\'er regime. In this work, and without loss of generality (by symmetry), we will consider only the case $\gamma>\gamma_c(h)$. To be precise we will rather consider the case $\gamma\ge\gamma_c(h)$ because the phenomenology observed 
for $\gamma>\gamma_c(h)$, that is for $M-\gamma_c(h) N \ge cN$ for some $c>0$ persists  also in a part of the window $M-\gamma_c(h) N =o(N)$ and we will analyze the model also in this window. %OK?
In different terms: the analysis in the Cram\'er regime is a Large Deviations analysis, but  the whole non-Cram\'er regime 
is equivalent from the Large Deviations viewpoint (the issues there are about sharp deviations). So there isn't much conceptual difference between $M-\gamma_c(h) N \ge cN$ and $M-\gamma_c(h) N =o (N)$, up to when 
$M-\gamma_c(h) N$ grows too slowly, as we shall see.

 \smallskip

Crucial for us is 
the probability distribution $\hat{\mathrm{K}}_h(\cdot, \cdot)$   defined by
\begin{equation}
\label{eq:hatK}
\hat{\mathrm{K}}_h(n,m) = K(n+m) e^{ h - n \tn(h) } \,,
\end{equation}
which, as announced informally just above,  allows to write the partition function as 
\begin{equation}
\label{eq:Zchat}
Z^c_{N,M,h}  =\exp (N \tn(h) ) \bP \left( (N,M) \in \hat{\tau}_h \right)\,,
\end{equation}
where $\hat{\tau}_h$ is the bivariate renewal process with inter-arrival distribution ${\hat{K}_h}(\cdot,\cdot)$. 
Next, we are going to have a closer look at this renewal process.

%The equation $\sum_{n,m} \hat{\mathrm{K}}_h(n,m) = 1$, that is the implicit equation for $\tn(h)$

\subsection{On the bivariate renewal $\hat{\tau}_h$}
\label{sec:brt}

 Let us write for conciseness  $\tn(h)=\tn_h$ (a practice that we will pick up again in the proofs), and drop the dependence on $h$ in $\hat \tau_h$:  $\hat\tau = (\hat\tau^{(1)}, \hat\tau^{(2)})$. In view of \eqref{eq:hatK}, it is clear that the distribution of this process is not symmetric, we have the marginals
\begin{equation}
\label{eq:marginals}
\begin{split}
\bP \Big( \htx = n \Big) \, &= \, \sum_{m \ge 1} K(n+m) \exp (h - \tn_h n) \, \stackrel{n \to \infty}\sim \,  \frac{\exp(h)}{(1+\ga)} \, \frac{L(n)}{n^{1+\ga}} \,  e^{-  \tn_h n} \,, \\
\bP \Big( \hty = m \Big) \, &= \, \sum_{n \ge 1} K(n+m) \exp (h - \tn_h n) \, \stackrel{m \to \infty}\sim \, \frac{\exp(h)}{\exp(\tn_h)-1} \, \, \frac{L(m)}{m^{2+\ga}}   \,.
\end{split}
\end{equation}
Let us also denote (the dependence in $h$ is implicit)
\begin{equation}
\hat \mu_1:=\bE[\htx_1] <+\infty , \qquad \hat\mu_2:=\bE[\hty_1]<+\infty\, , \end{equation}
so that $\gamma_c (h)= \hat\mu_2/\hat\mu_1$, cf. \eqref{eq:gammac}. 

\smallskip

We notice that the process $\htx$ has moments of all orders, and so $\{\htx_n\}_{n=0,1, \ldots}$ is in the domain of attraction of a normal law: we denote $a_n^{(1)}:=\sqrt{n}$ the scaling sequence for $\htx_n$. On the other hand, the process $\{\hty_n\}_{n=0,1, \ldots}$ is in the domain of attraction of an $\ga_2$-stable law, with $\ga_2:=(1+\ga)\wedge 2 >1$: its scaling sequence $a_n^{(2)}$ verifies
\begin{equation}
\label{def:an}
%\begin{split}
L(a_n^{(2)}) (a_n^{(2)})^{-\ga_2} \sim \,  \frac{1}{n} \  \text{ if } \ga_2<2
\quad \quad \text{ and } \quad \quad
\sigma(a_n^{(2)}) (a_n^{(2)})^{-2} \sim \, \frac 1n \ \text{ if } \ga_2=2
%\end{split}
\end{equation}
where 
\begin{equation}
\label{eq:sigma}
\sigma(n)\,:=\,\bE \Big[ \Big(\hty_1 \Big)^2 \ind_{\{\hty_1 \leq n\}} \Big]\,,
\end{equation}
 and diverges as a slowly varying function if $\bE[(\hty_1)^2]=+\infty$ (with $\sigma(n)/L(n) \to +\infty$ \cite{cf:RegVar}).
 In particular, $a_n^{(2)}$ is regularly varying with exponent $1/ \ga_2= (1+\ga)^{-1} \vee (1/2)$.

As an additional relevant definition, we  select a sequence $\{m_n^{(2)}\}_{n=1,2, \ldots}$ satisfying
\begin{equation}
\label{def:m_n}
\bP \big( \hty_1 > m_n^{(2)} \big) \stackrel{n\to\infty}\sim \frac1n ,
\end{equation}
so that $m_n^{(2)}$ gives the order of $\max_{1\leq j\leq n} \{\hty_j -\hty_{j-1}\}$. We stress that $m_n^{(2)}$ is regularly varying with exponent $(1+\ga)^{-1}$, and  that $m_n^{(2)}/ a_n^{(2)}\in [1/c,c]$ for some $c\ge 1$  if $\ga<1$, but $m_n^{(2)}/a_n^{(2)} \to 0$ when $\ga\geq 1$: in any case, there is a constant $c>0$ such that $m_n^{(2)} \leq c a_n^{(2)}$ for every $n$, i.e. $m_n^{(2)} =O( a_n^{(2)})$.

\medskip
Let us also stress that the bivariate renewal process $\hat{\tau}_h$ falls in the domain of attraction of an $(\ga_1=2,\ga_2)$ stable distribution (see e.g. \cite{cf:RG79} or \cite{cf:HOR84}).
We have, as $n\to\infty$,  that $\Big\{  \Big( \frac{\htx_n - \hat \mu_1 n}{ a_n^{(1)}}, \frac{\hty_n - \hat \mu_2 n}{a_n^{(2)}}  \Big)\Big\}_{n=1,2, \ldots}$ converges in distribution to $Z$, a non-degenerate $(2,\ga_2)$-bivariate stable law. Let us mention that in  \cite{cf:RG79} it is proven that:% $(\ga_1,\ga_2)$-stable distributions:
\begin{itemize}
\item[-] If $\ga_2 =2$ (i.e. $\ga \geq 1$), then $Z$ is a bivariate normal distribution.
\item[-] If $\ga_2 <2$ (i.e.  $\ga\in(0,1)$), then $Z$ is a couple of independent normal and $\ga_2$-stable distributions.
\end{itemize}

We mention that a bivariate local limit theorem is given in \cite{cf:DonLLT} and multivariate ($d$-dimensional) renewals are further studied in \cite{cf:B}: local large deviation estimates are given, as well as strong renewal theorems, i.e.\  asymptotics of $\bP( (n,m) \in \tau)$ as $(n,m)\to \infty$, when $(n,m)$ is \emph{close to} the \emph{favorite direction}  -- the favorite direction
exists when $\bE [\tau_1]$ is finite and it is 
the line $ t\mapsto t\bE[\tau_1]$, and \emph{close to} means at distance of the order of the  fluctuations around that direction -- we refer to \cite{cf:B} for further details (estimates when $(n,m)$ is \emph{away from} the favorite direction are also given).

\subsection{Non-Cram\'er regime and big-jump domain}
\label{sec:bigjump}
%Recall that $M-\gamma_c N \sim (\gamma-\gamma_c) N$, and that we do not want to study only the case $\gamma>\gamma_c$.

We drop the dependence of $\gamma_c(h)$ on $h$, 
and we set 
\begin{equation}
\label{def:tn}
t_N\, := M -\gamma_c N\, .
\end{equation}
Of course, having $\gamma>\gamma_c$ means that $t_N/N\ge c$ for some $c>0$. But it is natural and essentially not harder to 
tackle the problem 
assuming only
 %and we assume that $\{t_N\}_{n=1,2, \ldots}$ is regularly varying (NEEDED? I would take this away)
 \begin{equation}
 \label{hyp:bigjump1}
 t_N/a_N^{(2)}\to +\infty \ \text{ as } N\to\infty,
 \end{equation}
 with additionally, in the case $\ga\geq 1$ (recall the definition of $\sigma(n)$ after \eqref{def:an})
\begin{equation}
\label{hyp:bigjump}
\left( \frac{t_N}{a_N^{(2)}} \right)^2 \ \frac{\sigma(a_N)}{ \sigma(t_N)} \stackrel{N\to\infty}\sim \frac{t_N^2}{N\sigma(t_N)}  \geq C_0 \log N \ \text{ for a suitable choice of } C_0>0\, .
\end{equation}
If $\ga> 1$, as well as if $\ga=1$ and $\gs(n)=O(1)$ (i.e. if  $\bE[(\hty_1)^2]<\infty$), \eqref{hyp:bigjump} simply means  that $t_N \geq C'_0 \sqrt{N\log N}$ with $C'_0$ easily related to $C_0$. Note also that \eqref{hyp:bigjump1} implies $t_N \gg  \sqrt{N\log N}$ if
  $\bE[(\hty_1)^2]=\infty$.

We stress that the constants $C_0$ depends only on $K(\cdot)$ and, for the interested reader, it can be made explicit by tracking the constants in \eqref{eq:forC0_1} and \eqref{eq:forC0_2} where the value of $C_0$ is used.
This assumption is made to be sure that we lie in the so called \emph{big-jump domain}, as studied for example in the one-dimensional setting in~\cite{cf:DDS08}:
in our model it simply means that deviations -- and the event we focus on is $(N, M) \in \hat\tau$ --
are realized by an atypical deviation on just one of the increment variables $\hat\tau_{i+1}-\hat\tau_i$.
%, i.e.  by one single big jump and more precisely, just by one of the two components.
As we shall, this happens just under the assumption \eqref{hyp:bigjump1} for $\ga <1$ and this condition
is optimal (see Appendix \ref{app1}). For the case $\ga \ge 1$ the extra condition \eqref{hyp:bigjump} is not far from being optimal, but it is not: we discuss this point in Appendix \ref{app2}, but we do not treat it in full generality because it is a technically demanding issue that leads far from our main purposes.

\subsection{Mais results I: polymer trajectories}

We are now going to introduce two fundamental events in an informal, albeit  precise, fashion. The two events will be rephrased in a 
more formal way in \eqref{eq:BLBS}, once further notations will have been introduced. 
Choose sequences of positive numbers $\{u_N\}_{N=1,2, \ldots}$, $\{m^+_N\}_{N=1,2, \ldots}$, $\{a^+_N\}_{N=1,2, \ldots}$  and 
$\{\tilde a^+_N\}_{N=1,2, \ldots}$ such that
\begin{equation}
\label{eq:sequences}
u_N \gg 1\, ,
\quad 
t_N \gg 
m^+_N \gg m^{(2)}_N\, , 
  \quad
  t_N \gg a^+_N \gg a_N^{(2)}
\quad \text{ and } \quad
t_N \gg  \tilde a^+_N \gg a_N^{(2)}
\, .
\end{equation}
In practice, and to optimize the result that follows, $u_N$, $ m^+_N / m^{(2)}_N $, $a^+_N / a_N^{(2)}$ and 
$\tilde a^+_N / a_N^{(2)}$ should be chosen tending to $\infty$ in an arbitrarily slow fashion. 
%In view of the definition of $m_N$, see \eqref{def:m_n} and the observations right after \eqref{def:m_n}, we have also that $t_N \gg m_N=O(a_N^{(2)})$.

\smallskip

We then define the \emph{Big Loop} event $E^{(N)}_{\mathtt{BL}}$ to be the set of trajectories such that

\begin{enumerate}
\item there is one loop of size   larger than $t_N-a^+_N$ and smaller than $t_N+a^+_N$, so that, to leading order, it is of size $t_N$;
\item all other loops are smaller than $m^+_N$ (hence there is only one largest loop);
\item the length of neither of the two unbound strands is larger than $u_N$.
\end{enumerate}

\smallskip

The {\sl (large or macroscopic) Unbound Strand} event $E^{(N)}_{\mathtt{US}}$ is instead the set of trajectories such that

\begin{enumerate}
\item all loops are smaller than $m^+_N$;
\item the length of the unbound portion of the shorter strand does not exceed $u_N$;
\item the length of the unbound portion of the longer strand is larger than $t_N-\tilde a^+_N$ and smaller than $t_N+\tilde a^+_N$, so that, to leading order, it is of size $t_N$. 
\end{enumerate}

\begin{figure}[h] 
\centering
   \begin{minipage}[b]{.20\linewidth}
    \scalebox{0.4}{\includegraphics*{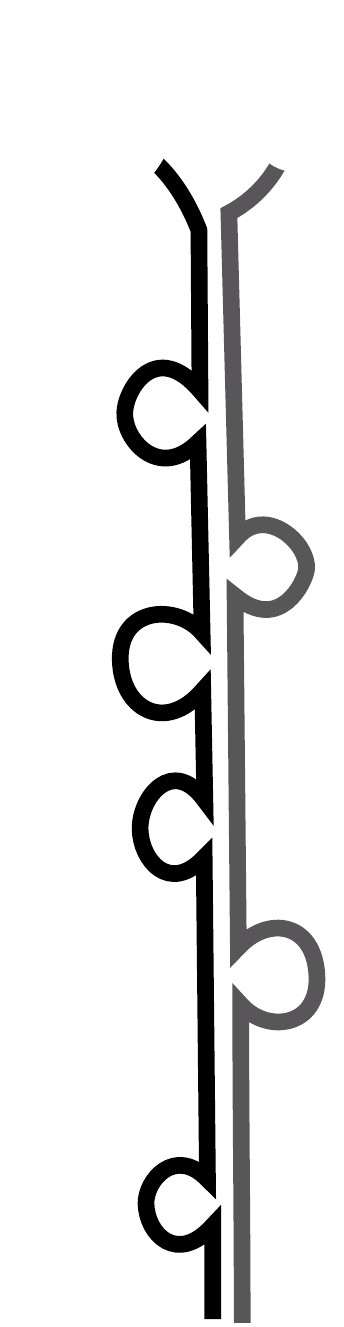}}
   \end{minipage}
     \quad
 \begin{minipage}[b]{.20\linewidth}
  \scalebox{1.3}{\includegraphics*{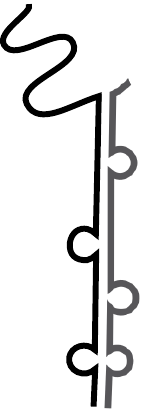}}
 \end{minipage}
   \quad
 \begin{minipage}[b]{.20\linewidth}
  \scalebox{1.3}{\includegraphics*{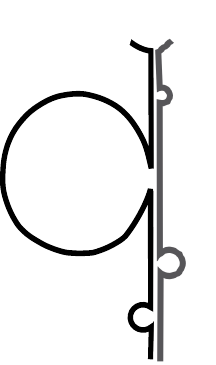}}
 \end{minipage}\hfill 
\caption{\label{fig:localized} Schematic image of the observable trajectories of the free gPS model in the Cram\'er regime (left), and in the non-Cram\'er regime (cf. Theorem~\ref{th:main-traj}): the Large Unbound Strand event (center, occuring when $\overline{\ga}  <\ga+1$) and the Big Loop event (right, occuring when $\overline{\ga}  > \ga +1$). In the constrained case the Unbound Strand event is not observed, and the free tails are of course absent.
What cannot be appreciated in this schematic view is the fact that the small loop distribution has exponential tail in the Cram\'er regime (hence the largest is $O(\log N)$) and that it has power law tail in the non-Cram\'er regime (hence the largest is $O(N^{a})$ for some $a\in (0,1)$: $O(m_N)$ to be precise).}
\end{figure}

\smallskip
Note that  $E^{(N)}_{\mathtt{BL}}\cap E^{(N)}_{\mathtt{US}}=\emptyset$ except, possibly, for finitely many $N$: the two conditions~(1) are incompatible. We refer to Fig.~\ref{fig:localized} for a schematic image of these two events.
\medskip

\begin{theorem}
\label{th:main-traj}
Under assumptions \eqref{hyp:bigjump1}, \eqref{hyp:bigjump} and \eqref{eq:sequences} we have that 
\begin{equation}
\label{eq:main-traj.1}
\lim_{N \to \infty} \bP^f_{N, M, h} \Big( E^{(N)}_{\mathtt{BL}}\cup E^{(N)}_{\mathtt{US}} \Big)\, =\, 1\,.
\end{equation}
Moreover
\begin{enumerate}
\item If $\overline{\ga}  <1$ (and hence $\sum_j K_f(j)=\infty$) then 
\begin{equation}
\label{eq:main-traj.2}
\lim_{N \to \infty} \bP^f_{N, M, h} \Big( E^{(N)}_{\mathtt{US}} \Big)\, =\, 1\,.
\end{equation}
\item If $\sum_j K_f(j)<\infty$ (and hence $\overline{\ga}  \ge 1$) then 
\begin{equation}
\label{eq:main-traj.3}
 \bP^f_{N, M, h} \Big(E^{(N)}_{\mathtt{US}}\Big)\stackrel{N \to \infty}\sim
 \frac{1}{1+Q_N}\ \ \text{ and} \ \
\bP^f_{N, M, h} \Big(E^{(N)}_{\mathtt{BL}}\Big)\stackrel{N \to \infty}\sim
 \frac{Q_N}{1+Q_N}\, ,
\end{equation}
with 
\begin{equation}
Q_N\, :=\, c_h N\, (t_N)^{-(2+\ga) + \overline{\ga}}\, \frac{L(t_N)}{\overline{L}(t_N)} \ \ \text{ and  } 
 c_h:= \frac{e^h  \sum_{j=0}^\infty K_f(j)}
 {\hat \mu_1(e^{\tn(h)}-1)} \, .
 \end{equation}
 \end{enumerate}
\end{theorem}

For conciseness the case $\overline{\ga} =1$ with $\sum_j K_f(j) = +\infty$ is not included in Theorem~\ref{th:main-traj}, but its is treated in full in Appendix~\ref{app:baralpha}. It is a marginal case in which an anomalous behavior appears:
 a big loop and a large unbound strand 
may coexist.

It is worth pointing out that, {\sl in most of the cases}, the expressions  in \eqref{eq:main-traj.3}  have a limit -- at least if $\{t_N\}_{N=1,2, \ldots}$  is not too wild (regularly varying is largely sufficient) -- and it is either one or zero.
In particular  when $t_N\sim c N$ for some $c>0$  we have
\begin{equation}
\lim_{N \to \infty}
 \bP^f_{N, M, h} \left(E^{(N)}_{\mathtt{US}}\right)\, =\, 
\begin{cases}
1 & \mathrm{if} \ \overline{\ga} <\ga+1\
\mathrm{or \ if} \ \overline{\ga} =\ga+1\ \mathrm{and} \ \overline{L}(t) \stackrel{t \to \infty} \gg L(t)  ,\\ 
0 & \mathrm{if} \ 
 \overline{\ga} >\ga+1\ \mathrm{or \ if} \
\overline{\ga} =\ga+1\ \mathrm{and} \ \overline{L}(t) \stackrel{t \to \infty} \ll L(t)  . \\ 
\end{cases}
\end{equation}
(This is true also in the case $\overline{\ga}=1$ with $\sum_j K_f(j) =+\infty$, see \eqref{eq:tildeQ}-\eqref{eq:Palphabar}).
Note that in the case in which $\overline{\ga} =\ga+1$ and the ratio of the two slowly varying function has a limit which is neither $0$ nor $\infty$, the limit of the probability of the unbound strand event exists and it is an explicit value in $(0,1)$. 
%This shows also that cases in which the limit does not exist are easily built, but they are very special cases because they come up only when the slowly varying functions play a role. 

\section{Main results II: sharp estimates on the partition functions}

In this section, we give the asymptotic behavior  of $e^{-N \tn_h} Z_{N,M,h}^c =\bP( (N,M) \in \hat\tau)$  in the big-jump domain.
Then we present the asymptotic behavior  of $Z_{N,M,h}^f$. Both in the constrained and free case we also give more  technical estimates that identify some events to which we can restrict the partition functions  without modifying them in a relevant way. 
Theorem~\ref{th:main-traj} turns out to be a corollary of these technical estimates, as we explain in the final part of the section.
 
\smallskip
In this section and in the rest of the paper we deal with order statistics and we introduce here the relative definitions.
Consider the (non-increasing) order statistics $\big\{\cM_{1,k}, \cM_{2,k}, \ldots, \cM_{k,k}\big\}$ of the IID family
$\{ \hty_j-\hty_{j-1}\}_{j=1, \ldots, k}$. In particular $\cM_{1,k}$ is a maximum of this finite sequence.  
We will consider the order statistics also for $k$ random, notably for $k= \kappa_N:=\max\{ i\,:\, \htx_i \in [0, N]\}$.

\subsection{On the constrained partition function}

We start with an important estimate for the constrained partition function (more precisely for the renewal mass function $\bP((N,M)\in\hat \tau)$), that is essential for the study of the free partition function, as one can imagine from its definition \eqref{eq:zf}. 
It is worth insisting on the link between $\bP$ and the measure we are interested in for the constrained case: 
\begin{equation}
\label{eq:PcP}
\bP_{N, M, h}^c(\cdot)\, =\, 
\bP( \, \cdot \, \vert (N, M) \in \hat \tau )\, .
\end{equation}

\begin{theorem}
\label{th:aimZc}
Assume that $\ga >0$ and \eqref{hyp:bigjump1}. Moreover if $\ga\geq 1$ assume also \eqref{hyp:bigjump}.
Then (recall that $M=\gamma_c N+t_N$) we have that
\begin{equation}
\bP \big(  (N,M) \in \hat\tau \big) = \bP \Big(\htx_{\kappa_N}=N, \, \hty_{\kappa_N}= M \Big) \stackrel{N\to\infty}{\sim} \frac{N}{\hat \mu_1^2} \bP \Big(\hty_1 = \lceil t_N  \rceil \Big).
\label{eq:hatp}
\end{equation}
Moreover, for every $\eta \in (0,1)$ there exist $\gep_0>0$  such that for every $\gep\in(0, \gep_0)$ and  $N$ sufficiently large (how large may depend on $\gep$), we have 
\begin{multline}
\label{eq:aimlb}
\bP\bigg( 
\cM_{1,\kappa_N} - t_N  \in \Big[-\frac{a_N^{(2)}}\gep , \frac{a_N^{(2)}}{\gep} \Big], \, 
\cM_{2,\kappa_N} \le \frac {m_N^{(2)}}{ \gep}, \, \htx_{\kappa_N}=N, \, \hty_{\kappa_N}= M
\bigg)\, \ge \\  \ge  (1-\eta) \frac{N}{\hat \mu_1^2} \bP \Big( \hty_1= \lceil t_N  \rceil\Big)\,.
\end{multline}
\end{theorem}

\subsection{On the free partition function}

We now give the behavior of the free partition function and identify trajectories contributing the most to it.
Let us introduce some notations:
\begin{equation}
V_1^{(N)} := N - \htx_{\kappa_N} ,\qquad
V_2^{(N)} := M - \hty_{\kappa_N} \, ,
\end{equation}
the  lengths of the free parts, see Fig~\ref{fig:explain}.

For a set $A$ of allowed trajectories, we define $Z_{N,M,h}^f(A)$ the partition function restricted to trajectories in $A$ (by restricting the summation over subsets of $\{1, \ldots N\}\times  \{1, \ldots M\}$ to those in $A$). For example, $Z_{N,M,h}^f( (N,M)\in \hat \tau) = Z_{N,M,h}^c$.

 We set $\overline{K}:= \sum_{j=1}^\infty K_f(j)$ if the sum is finite, and we set $\overline{K}=0$ if $\sum_j K_f(j) = +\infty$. 

%Since $Z_{N,M, h}^{f}$ is obtained by summing over subsets of $\{1, \ldots N\}\times  \{1, \ldots M\}$
%that coincide with the trajectories of $\hat \tau$, the notion of the partition function restricted to any set of allowed trajectories of
%$\hat \tau$ in $\{1, \ldots N\}\times  \{1, \ldots M\}$ is well defined.

\begin{theorem}
\label{th:aimZf}
Assume that $\ga >0$ and \eqref{hyp:bigjump1}. Moreover  if $\ga\geq 1$ assume also \eqref{hyp:bigjump} and 
if  $\overline{\ga} =1$ assume  that $\sum_j K_f(j) <\infty$.
Then for $N \to \infty$
\begin{multline}
\label{eq:aimZf1}
e^{- N \tn(h) }  Z_{N,M, h}^{f} = (1+o(1)) 
\overline{K} \,
\frac{N}{\hat \mu_1^2 } \bigg(\sum_{i\ge 0} K_f(i) e^{- i \tn(h)} \bigg) 
\bP \Big( \hty_1= \lceil t_N \rceil \Big)\\
 + (1+o(1)) \frac{1}{ \hat \mu_1}\bigg(\sum_{i\geq 0} K_f(i) e^{- i \tn_h} \bigg)  K_f\left( t_N\right)\, .
\end{multline}
Moreover, for every $\eta \in (0,1)$ there exist $\gep_0>0$  such that for every $\gep\in(0, \gep_0)$ and  $N$ sufficiently large (how large may depend on $\gep$), such that
\begin{align}
e^{- N \tn_h } Z_{N,M , h}^f \bigg( V_1^{(N)}, V_2^{(N)} \leq \frac1\gep ,\, & \cM_{1,\kappa_N} \in \Big[ t_N -\frac{a_N^{(2)}}\gep , t_N + \frac{a_N^{(2)}}\gep  \Big] , \, 
\cM_{2,\kappa_N} \le \frac {m_N^{(2)}}{ \gep} \bigg) 
\notag\\
& \geq  (1-\eta)\overline{K}\, 
 \frac{N}{\hat \mu_1^2 }\bigg(\sum_{i\geq 0} K_f(i) e^{- i \tn_h } \bigg) \bP\left( \hty_1=t_N \right) \, ;
 \label{eq:aimZf2.1}
 \\
e^{- N \tn_h } Z_{N,M , h}^f \bigg(V_1^{(N)}\leq \frac1\gep , V_2^{(N)} &\in \Big[ t_N -\frac{a_N^{(2)}}{\gep}, t_N+ \frac{a_N^{(2)}}{\gep} \Big], \cM_{1,\kappa_N}\leq \frac{m_N^{(2)}}{\gep} \bigg) \notag
\\
& \geq  (1-\eta)\frac{1}{\hat \mu_1}\bigg(\sum_{i\geq 0} K_f(i) e^{- i \tn_h} \bigg)  K_f\left( t_N\right)\, .
\label{eq:aimZf2.2}
\end{align}
\end{theorem}

We remark that when $\overline{\ga} <1$, that is $\overline{K} =0$, the right-hand side of
\eqref{eq:aimZf1} reduces to one term and \eqref{eq:aimZf2.1} becomes trivial. 
The case $\overline{\ga}  =1$ with $\sum_j K_f(j)=+\infty$ is treated in Theorem~\ref{th:baralpha}.

%
% In fact, the asymptotics of \eqref{eq:aimZf1} is modified: the first term $\overline K$ has to be replaced by $\overline K(t_N):= \sum_{j=1}^{t_N} K_f(j)$ -- which diverges to $+\infty$ as a slowly varying function. Moreover, the main contribution of this term does not come from \eqref{eq:aimZf2.1}, but from a set of trajectories with both a big loop and a large unbound strand.
%The details of this case are developed in Appendix \ref{app:baralpha}.

\subsection{Back to the Big Loop and Unbound Strand events. %the proof of Theorem~\ref{th:main-traj}
} 
\label{sec:BLUS}

The notations we have introduced
 allow a  compact formulation of the two  key events of  Theorem~\ref{th:main-traj}:

\begin{equation}
\begin{split}
\label{eq:BLBS}
E^{(N)}_{\mathtt{BL}}\, &=\, 
\left\{
\cM_{1,\kappa_N}\in [t_N-a_N^+, t_N+a_N^+], \, 
\cM_{2,\kappa_N}< m^+_N, \, 
\max\left(V_1^{(N)}, V_2^{(N)}\right) \le u_N
\right\}\, ,
\\
E^{(N)}_{\mathtt{US}}\, &=\, 
\left\{
\cM_{1,\kappa_N}< m^+_N, \, 
V_1^{(N)} \le u_N,
\, 
V_2^{(N)} \in [t_N-\tilde a_N^+, t_N+\tilde a_N^+]
\right\}\, .
\end{split}
\end{equation}

\noindent
\emph{Proof of Theorem~\ref{th:main-traj}.}
This is just a book-keeping exercise using the three estimates in Theorem~\ref{th:aimZf}, together with the definition of $K_f(t_N)$ and the estimate of $\bP(\hty_1 = \lceil t_N \rceil)$ in \eqref{eq:marginals}. 
\qed

\smallskip

Finally, in the same way, we obtain from Theorem~\ref{th:aimZc} 
the following complement to Theorem~\ref{th:main-traj} for the constrained case.

\begin{theorem}
\label{th:main-traj-c}
Under assumptions \eqref{hyp:bigjump1} (and additionally \eqref{hyp:bigjump} if $\ga\ge 1$) and \eqref{eq:sequences} we have that 
\begin{equation}
\label{eq:main-traj-c}
\lim_{N \to \infty} \bP^c_{N, M, h} \left(E^{(N)}_{\mathtt{BL},0}\right)\, =\, 1\,,
\end{equation}
where 
$E^{(N)}_{\mathtt{BL},0}$ is the event
$E^{(N)}_{\mathtt{BL}}$ with the more stringent condition that $\max (V_1^{(N)},V_2^{(N)})=0$.
\end{theorem}

\subsection{A word about the arguments of  proof and organization of the remaining sections}

As we pointed out at the beginning of the introduction, \emph{condensation phenomena} are widely studied in the mathematical literature
(see \cite{cf:DDS08,cf:AL} and references therein), but not in the multivariate context. The full multivariate context is the object of \cite{cf:B}, where renewal estimates $\bP((n_1,\ldots, n_d)\in \hat \tau)$ are given: in particular, in the big-jump domain, only rough (general) bounds are given. Here, we address only the special bivariate case motivated by the application so that we are able to give the exact asymptotic behavior.
One of the main 
 difficulties we face is that, on the event $(N,M)\in\hat \tau$, the number $\kappa_N$ of renewal points is random and highly constrained by this event. We show that in the big-jump domain considered in Section \ref{sec:bigjump}, the main contribution to the probability $\bP((N,M)\in\hat \tau)$ comes from trajectories with a number of renewal points that is approximately $k_N = N/\hat\mu_1 +O( \sqrt{N}) $. For this number $k_N$,  $\hat\tau^{(1)}_{k_N}$ does not have to deviate from its typical behavior to be equal to $N$, but $\hat\tau_{k_N}^{(2)}$ has to deviate from its typical behavior
 to reach $M$ 
 and it does so
  by making one single big jump, of order $t_N +O( a_N^{(2)})$. In this sense, if we accept that $\kappa_N$ is forced to be $N/\hat \mu_1 +O( \sqrt{N})$ by the condition $N\in \hat \tau^{(1)}$, we can focus on $M\in \hat\tau^{(2)}$ and the problem becomes 
 \emph{almost one dimensional}. This turns out to be  a lower bound strategy: for a corresponding upper bound 
 we have to show that all other trajectories bring a negligible contribution  to $\bP((N,M)\in \hat\tau)$.

\smallskip

In the rest of the paper, we estimate separately the constrained and free partition functions. We deal with the constrained partition function in Section \ref{sec:constrained}:
the main term \eqref{eq:aimlb} in Section~\ref{sec:constrained1} and the remaining negligible contributions in Section \ref{sec:constrained2}.
The free partition function is dealt with in Section~\ref{sec:free}: the main terms \eqref{eq:aimZf2.1} and \eqref{eq:aimZf2.2} in Section~\ref{sec:free1} and the remaining negligible contributions in Section \ref{sec:free2}.
In Appendix \ref{app:baralpha} we complete the analysis of the case $\overline{\ga}  =1$.
In Appendix \ref{app} we discuss the transition from the big-jump regime (a single big jump, with a big deviation of just one of the two components) to 
 the  Cram\'er deviation strategy (no big jump). %and we claim, but do not prove, that for $\ga \ge 1$,  there is a tiny transition region in which  there is one single big jump but performed by both components.

\smallskip

To keep things simpler in the rest of the paper, and with some abuse of notation, we will systematically omit the integer part in the formulas.

\section{The constrained partition function: proof of Theorem~\ref{th:aimZc}}
\label{sec:constrained}

\subsection{Proof of the lower bound \eqref{eq:aimlb}}
\label{sec:constrained1}

We start by decomposing the event of interest according to  $\kappa_N=k$. The probability of such an event, restricted to $\{\kappa_N=k\}$, becomes (recall that $M=\gamma_c N +t_N$)
\begin{multline}
\label{eq:tildeP1}
\bP \Big( 
\cM_{1,k} \ge t_N -\frac{a_N^{(2)}}\gep , \, 
\cM_{2,k} \le \frac {m_N^{(2)}}{ \gep}, \, \htx_k=N, \, \hty_k= M
\Big)
\, =
\\
\bP \bigg( \bigcup_{j=1}^k
\Big \{ \hty_j-\hty_{j-1} \in t_N+ I_N , \, 
\max_{i \neq j}\left( \hty_i- \hty_{i-1}\right) \le
 \frac {m_N^{(2)}}{ \gep}, \, \htx_k=N, \, \hty_k= M
\Big \}
\bigg)\, ,
\end{multline}
where we defined $I_N :=\{-a_N^{(2)}/\gep, \ldots, a_N^{(2)}/\gep\}$.
Since $t_N -a_N^{(2)}/\gep$ is larger than 
${m_N^{(2)}}/{ \gep}$ (recall that $m_n^{(2)}\leq c a_n^{(2)}$ for every $n$  and \eqref{hyp:bigjump}) for $N$ sufficiently large, the union in the right-hand side of \eqref{eq:tildeP1}
is a union of disjoint events that have all the same probability. This term is equal to
\begin{multline}
\label{eq:union}
k \, \bP \Big( \hty_1 \in t_N+I_N , \, 
\max_{i=2, \ldots, k}\left( \hty_i-\hty_{i-1}\right) \le 
 \frac {m_N^{(2)}}{ \gep}, \, \htx_k=N, \, \hty_k= \gamma_c N +t_N
\Big) \, = \\
k \, \sum_{y \in I_N } \sum_{x\in \bbN} \bP \Big(\hty_1 =  t_N  +y, \, \htx_1=x
\Big) 
\bP \Big(\cM_{1,k-1} \leq  \frac {m_N^{(2)}}{ \gep},  \htx_{k-1}=N-x, \, \hty_{k-1}=  \gamma_c N -y
\Big)\, ,
\end{multline}
where we have used that $\{(\htx_j-\htx_1,\hty_j-\hty_1)\}_{j=2, \ldots,k }$ and  
$\{(\htx_j,\hty_j)\}_{j=1, \ldots, k-1 }$ have the same law.
% Moreover we have used that $M- \lceil t_N\rceil= \lfloor \gamma_c N\rfloor$.

Since we are after a lower bound we may and do restrict the sum over $x$ between $1$ and $1/\gep$
and $y \in I_N:=\{-a_N^{(2)}/\gep, \ldots, a_N^{(2)}/\gep\}$. And using that $\bP(\hty_1 =n) $ is regularly varying, we have that uniformly for such $x$ and $y \in I_N$
\[
 \bP \Big(\hty_1 =  t_N +y, \, \htx_1=x\Big) \, \ge 
 \left(1- \gd_N \right)  \bP \Big(\hty_1 =  t_N 
\Big)  \tilde \bP \Big(\htx_1=x \, \Big \vert \, \hty_1 =   t_N  +y
\Big)\, ,\]
where $\gd_N =\gd_N(\gep)\ge 0$ is such that $\lim_{N \to \infty}\gd_N=0$. 
If now we set $p_h(x):= (e^{\tn_h}-1) e^{- x\tn_h}$, by using \eqref{eq:marginals} we have that 
\begin{equation}
\begin{split}
\bP \Big(\htx_1=x \, \big \vert \, \hty_1 =  &  t_N  +y
\Big)\, =\, \frac{\bP \big(\htx_1=x , \, \hty_1 =  t_N  +y
\big)}{\bP \big( \hty_1 =   t_N +y
\big)}
\\
& \ge \, \left(1-\gd_N\right)
\frac{p_h(x)K\left( t_N +x+y\right)}{K\left(  t_N +y\right)}\,
\ge \, \left(1-\gd_N\right)^2
p_h(x)
\, 
,
\end{split}
\end{equation}
%and $K(\lceil t_N\rceil+x+y)/K(\lceil t_N\rceil+y)$ is also bounded below by $1-\gd_N$, 
possibly for a different choice of $\gd_N =\gd_N(\gep)$. 
Therefore, going back to \eqref{eq:tildeP1} we see that (again, by redefining $\gd_N $)
\begin{multline}
\label{eq:tildeP2}
\bP \Big( 
\cM_{1,k} \in t_N + I_N, \, 
\cM_{2,k}\le \frac {m_N^{(2)}}{ \gep}, \, \htx_k=N, \, \hty_k= M
\Big)
\, \ge 
\\
(1-\gd_N) k \, \bP \big(\hty_1 =   t_N 
\big) 
 \sum_{x=1}^{1/\gep} p_h(x) 
\bP \Big(
\cM_{1,k-1} \leq  \frac {m_N^{(2)}}{ \gep},  \htx_{k-1}=N-x, \, \hty_{k-1}- \gamma_c N \in I_{N}
\Big)
\, .
\end{multline}
We now sum over the values of $k$ and we restrict to  $k \in  [(N/ \hat\mu_1)- \sqrt{N/\gep} ,
(N/ \hat\mu_1) +\sqrt{N/\gep}]\cap \bbZ:= J_N$. 
%For \eqref{eq:tildeP2} we have also used that $\tilde \bP( Y_1=\cdot)$ is regularly varying.
Hence,  redefining $\gd_N$, 
 the left-hand side of \eqref{eq:aimlb} is bounded from below by 
 \begin{equation}
 \label{eq:tildeP2.2.1}
 (1-\gd_N)  \bP \Big(\hty_1 =   t_N 
\Big) \frac N{\hat \mu_1} \sum_{x=1}^{1/\gep} p_h(x) P_\gep(x)\,  ,
 \end{equation}
 where we defined, with $n_N^+ := \max J_N$,
\begin{equation}
 P_{\gep}(x):= \sum_{k\in J_N}
\bP \bigg(
\max_{i=1, \ldots, n^+_{N}} \left( \hty_i-\hty_{i-1}\right)\leq 
 \frac {m_N^{(2)}}{ \gep},  \htx_{k-1}=N-x, \, \hty_{k-1}-  \gamma_c N \in I_{N}
\bigg)
\end{equation}
For $P_{\gep}(x)$, we observe right away that by introducing also
%$n^+_{N}:= \max  J_N$  and 
$n^-_{N}:= \min  J_N$ -- note that $n_N^{\pm}$ are equal to
 $(N/\hat \mu_1) \pm \sqrt{N/\gep}$ -- we have
\begin{equation}
\label{eq:tildeP2.2}
\begin{split}
P_\gep(x) \,  \ge \, &  
 \sum_{k\in J_N}
\bP \bigg( \cM_{1,n^+_N} \leq   
 \frac {m_N^{(2)}}{ \gep},  \htx_{k-1}=N-x
\bigg)
\\ 
&-
 \sum_{k\in J_N}
\bP \bigg( \cM_{1,n^+_N} \leq
 \frac {m_N^{(2)}}{ \gep},  \htx_{k-1}=N-x, \, \hty_{n^+_{N}}- \gamma_c N > \frac{a_N^{(2)}}{ \gep}
\bigg)
\\
 &-\sum_{k\in J_N}
\bP \bigg( \cM_{1,n^+_N} \leq  \frac {m_N^{(2)}}{ \gep},  \htx_{k-1}=N-x, \, \hty_{n^-_{N}}-  \gamma_c N < \frac{a_N^{(2)}}{ \gep}
\bigg)
\\ &\geq  \bP \left( E_1\cap E_2(x)  \right) - \bP \left( E^+_3 \right)- \bP \left( E^-_3 \right)
\, ,
\end{split}
\end{equation}
where 
\begin{equation}
E_1\, :=\, \Big \{ \cM_{1,n^+_N} \leq
 \frac {m_N^{(2)}}{ \gep}  \Big\}\, , \qquad 
 E_2(x)\, :=\, \left\{ \exists k \in  J_N \textrm{ such that } \htx_{k-1}=N-x\right\}\, ,
\end{equation}
and 
\begin{equation}
E_3^+\, :=\, \Big \{ \hty_{n^+_{N}}-  \gamma_c N  > \frac{a_N^{(2)}}{ \gep} \Big \}
\, , \qquad   E_3^-\, :=\, \Big \{ 
\hty_{n^-_{N}}-  \gamma_c N  < \frac{a_N^{(2)}}{ \gep}
\Big \}\, .
\end{equation}
We now estimate separately the probability of these events.

\subsubsection{$E_1$ has probability close to one}
For this, we use that $\bP( \hty_1 >n)$ is regularly varying with index $(1+\ga)^{-1}$ together with the definition \eqref{def:m_n} of $m_n^{(2)}$ to obtain that for $N$ larger than some constant $N_0=N_0(\gep)$ we have $\bP \left(\hty_1 > \tfrac1\gep m_N^{(2)}\right) \leq 2 \gep^{1+\ga} N^{-1}$.
Therefore, we have for $N\geq N_0$
\begin{equation}
\label{eq:est-E_1}
\bP \left( E_1\right) \, \geq \, \left( 1- 2 \gep^{1+\ga} N^{-1} \right)^{n^+_{N}} \geq e^{- 3 \gep^{1+\ga}} \, .
\end{equation}
where we used that $n_N^+ \leq N$ and $\gep$ small.

\subsubsection{$E_2(x)$ has probability close to $1/\hat \mu_1$}
The probability of $E_2(x)$ is estimated by writing
\begin{equation}
\label{eq:event2}
\begin{split}
\bP \left( E_2(x)\right)\, &= \bP(N-x \in \htx)  - \sum_{k < n_N^{-}} \bP \left( \htx_{k} = N-x \right) - \sum_{k >  n_N^+}  \bP \left( \htx_{k} = N-x \right) \\
&\geq \bP(N-x \in \htx)  -  \bP\Big( \htx_{n_N^{-}} \geq N -\tfrac1\gep \Big) - \bP \Big( \htx_{n_N^+} \leq N \Big)\, ,
\end{split}
\end{equation}
where for the second term we used that $\bP \big(\exists\, k< n_{N}^-  \text{ s.t. } \htx_{k} = N-x \big) \leq \bP \big( \htx_{n_N^-} \geq N-x  \big) $ together with the fact that $x \leq 1/ \gep$ (and similarly for the last term).

First, because the inter-arrivals of $\htx$ are exponentially integrable, $\vert  \bP(N \in \htx) - 1/\hat \mu_1\vert \le \exp(-c N)$
for $N\ge N_0$ with $c>0$ and $N_0$ that depend on the inter-arrival law \cite{cf:Kendall}. Therefore, uniformly in $x=1, \ldots,1/\gep$,
we have that for $N$ sufficiently large $\bP(N-x \in \htx) \geq 1/\hat \mu_1 - e^{- c N/2 }$.

For the remaining terms in \eqref{eq:event2} it is just a matter of using the Central Limit Theorem. In fact,
recalling that $n_N^- = N/\hat \mu_1 - \sqrt{N/\gep}$, we have
\begin{equation}
\label{eq:lowertau1}
\bP\Big(   \htx_{n_N^{-}} \geq N -\tfrac1\gep\Big) = \bP\Big(   \htx_{n_N^{-}} - \hat \mu_1 n_N^{-} \geq \hat \mu_1 \gep^{-1/2} \sqrt{N} - \tfrac1\gep \Big)
\leq e^{-c  \gep^{-1}},
\end{equation}
%which is a standard application of the Central Limit Theorem, 
for $N$ larger than some $N_0=N_0 (\gep)$.
On the other hand, we also have that
\begin{equation}
\label{eq:uppertau1}
\bP\Big(   \htx_{n_N^{+}} \leq N \Big) = \bP\Big(   \htx_{n_N^{+}} - \hat \mu_1 n_N^{+} \leq -  \hat \mu_1 \gep^{-1/2} \sqrt{N} \Big)
\leq e^{-c' \gep^{-1}},
\end{equation}
provided again that $N$ is large enough.

Therefore we have proven  that for every $\eta \in (0,1)$ there exists $\gep_0$ and $N_0:(0,1) \to \bbN$ such that for  $\gep \in (0, \gep_0)$ and $N \ge  N_0(\gep)$ we have
\begin{equation}
\label{eq:est-E_2}
\min_{x=1, \ldots, 1/ \gep}
\bP \left( E_2(x)\right)\, \ge\, \frac{1-\eta }{\hat \mu_1}\, .
\end{equation}

\subsubsection{$E_3^{\pm}$ have a small probability }
This is a consequence of the convergence to stable limit law.
In fact, using that $\gamma_c = \hat \mu_2/\hat \mu_1$ so that $\gamma_c N = \hat \mu_2 n_N^+ -\hat\mu_2 \sqrt{N/\gep} $, we get
\begin{equation}
E_3^+\, = \bigg\{ \hty_{n^+_N}- \hat \mu_2 n^+_N\, >\, \frac{a_N^{(2)}}{\gep} -  \hat\mu_2 \sqrt{N / \gep}\bigg\}
\subset \bigg\{ \hty_{n^+_N}- \hat \mu_2 n^+_N\, >\, \frac{a_N^{(2)}}{2\gep}\bigg\}\, ,
\end{equation}
where the last inclusion holds provided that $\gep$
is sufficiently small, since there is a constant $c$ such that $a_N^{(2)}\geq c \sqrt{N}$ for all $N$ (we actually simply need $N$ to be large if $a_N^{(2)}/\sqrt{N}\to +\infty$ for $N \to \infty$, which is the case when $\bE[(\hty_1)^2]=+\infty$). Very much in the same way we get also to
\begin{equation}
E_3^- \, \subset\, 
\bigg\{ \hty_{n^-_N}- \hat\mu_2 n^-_N\, <\, -\frac{a_N^{(2)}}{2\gep}\bigg\}\, .
\end{equation}
Since 
$(\hty_{n^\pm_N}- \hat \mu_2 n^\pm_N)/a_{n^\pm_N}^{(2)}$ converges in law for $N \to\infty$ to a stable limit variable $Y$, and using that  $a_N^{(2)}/a^{(2)}_{n_N^{\pm}} \to \hat \mu_1^{1/\ga_2}$ (since $n_N^\pm \sim N/ \mu_1$ and $a_N^{(2)}$ is regularly varying with exponent $\ga_2^{-1}$, recall $\ga_2:=\min(1+\ga,2)$), it is straightforward to see that
\begin{equation}
\label{eq:est-E_3}
\limsup_{N \to \infty } \bP(E_3^+)\, \le \, \bP \bigg( Y \ge \frac {\hat \mu_1^{1/\ga_2}}{2\gep} \bigg)
\ \text{ and } \
\limsup_{N \to \infty } \bP(E_3^-)\, \le \, \bP \bigg( Y \le- \frac {\hat \mu_1^{1/\ga_2}}{2\gep} \bigg)\, ,
\end{equation}
which are both vanishing as $\gep \searrow 0$. 

\medskip
We therefore see that
\eqref{eq:est-E_1}, \eqref{eq:est-E_2} and \eqref{eq:est-E_3} yield that, provided that $\gep_0$ is small enough, for every $\gep<\gep_0$ and $N$ large enough (how large depends on $\gep$), $ P_{\gep}(x) \geq (1-\eta)/\mu_1$ uniformly for $x\in\{1,\ldots, 1/\gep\}$.
If we now go back to \eqref{eq:tildeP2.2.1} and \eqref{eq:tildeP2.2}, and using that $\sum_{x\geq 1} p_h(x) =1$, we obtain\eqref{eq:aimlb}.
\qed

\subsection{Proof of \eqref{eq:hatp}}
\label{sec:constrained2}

In view of \eqref{eq:aimlb}, we simply need to give an upper bound on the probability $\hP\left( (N,\gamma_c N +t_N) \in\hat \tau \right)$. Fix some $\gep>0$.

\smallskip
\noindent
{\bf First step.}
We control
\begin{align}
\bP&  \Big( \cM_{1,\kappa_N} > (1-\gep)  t_N , \htx_{\kappa_N} = N, \hty_{\kappa_N} =\gamma_c N +t_N\Big) \notag\\
&  \le \sum_{k=1}^{N} k \sum_{y > - \gep t_N} \sum_{x\in\bbN}  \bP\left(\hty_1=   t_N  +y ,\, \htx_1 =x \right)  \bP\left( \htx_{k-1} = N-x , \hty_{k-1} =  \gamma_c N  - y \right)\, .
\label{eq:part1}
\end{align}

Recalling \eqref{eq:hatK} and \eqref{eq:marginals}, we have that there is some $N_0=N_0(\gep)$ and some $\eta=\eta_{\gep}$ (with $\eta_{\gep} \to 0$ as $\gep\downarrow 0$), such that for all $N\geq N_0$,  $x\leq 1/\gep$ and $y > - \gep t_N$ we have
\begin{equation}
\bP \Big(\hty_1=  t_N +y ,\, \htx_1 =x \Big) = \frac{L\left(  t_N +y+\right)}{\left(  t_N +y+x\right)^{2+\ga}} e^{h-x \tn_h}  \leq (1+\eta) \hP \Big(\hty_1= t_N \Big) p_h(x) ,
\end{equation}
where we recall that $p_h(x):= (e^{\tn_h}-1) e^{-x \tn_h}$.
Note that we also have that there is a constant $c$ such that uniformly for $x\in\bbN$
\begin{equation}
\label{eq:uniformbound}
\bP \Big(\hty_1= z ,\, \htx_1 =x \Big) \leq c L(z) z^{-(2+\ga)} p_h(x)\, .
\end{equation}
We can use that to get that uniformly for $y\geq -t_N/2$ (so that $t_N +y \geq t_N/2$) we have that for any $x\geq 1$
\begin{equation}
\bP \Big(\hty_1=  t_N +y ,\, \htx_1 =x \Big)\leq  c' p_h(x) \bP \Big( \hty_1 =   t_N  \Big).
\end{equation}

Then, dividing \eqref{eq:part1} according to whether $x\le 1/\gep$ or $x>1/\gep$ (and summing over $y > \gep t_N$), we obtain the following upper bound
\begin{multline}
\label{eq:sum<1/gep}
(1+\eta) \sum_{x=1}^{1/\gep} p_h(x) \bP \Big( \hty_1 =  t_N  \Big)  \sum_{k=1}^{N} k  \hP \Big( \htx_{k-1} = N-x \Big)\\
 + c \sum_{x=1/\gep}^{N} p_h(x) \bP \Big( \hty_1 =   t_N \Big)  \sum_{k=1}^{N} k \hP \Big( \htx_{k-1} = N-x \Big)\, .
\end{multline}
The second term is bounded from above by 
\begin{equation}
c N \hP \Big(\hty_1=  t_N \Big) \sum_{x>1/\gep} p_h(x) = c e^{-\tn_h/\gep } \times N \hP \Big(\hty_1=  t_N \Big)\, .
\end{equation}

In the first term \eqref{eq:sum<1/gep}, we split the sum according to whether $k$ is smaller or greater than $(1+\gep)N/\hat \mu_1$: we get that
\begin{equation}
\begin{split}
\sum_{k=1}^{N} k  \hP&\Big( \htx_{k-1} = N-x \Big)\\
& \le (1+\gep) \frac{N}{\hat \mu_1} \hP \Big( N- x \in \htx \Big) + N \hP \Big( \exists k > (1+\gep)N/\hat \mu_1 \ \text{ s.t. } \htx_{k-1} = N -x  \Big) \\
& \le (1+\gep)^2  \frac{N}{\hat \mu_1^2} + N \hP \Big( \htx_{(1+\gep)N/\hat \mu_1} \leq N -x \Big) \, ,
\end{split}
\end{equation}
where we used that in the first part $k \leq (1+\gep)N/\hat \mu_1$, and the renewal theorem to get that $ \hP \left( N-x \in\htx \right)\leq (1+\gep)N/\hat \mu_1$ uniformly for $x\le 1/\gep$ and $N$ large enough (how large depends on $\gep$).
The second term is exponentially small since it is a large deviation for $ \htx$ ($x$ here is bounded by  $1/\gep$).
Recalling that $\sum p_h(x) =1$, the first term \eqref{eq:sum<1/gep} is therefore bounded from above by
\begin{equation}
(1+\eta)\left( (1+\gep)^2+ e^{- c_\gep N} \right) \frac{N}{\hat \mu_1^2} \bP \left( \hty_1 =  t_N  \right)   \, ,
\end{equation}
In the end, the left-hand side of \eqref{eq:part1} is bounded by
\begin{equation}
\left(1+ \eta'_{\gep} \right)   \frac{N}{\hat \mu_1^2} \bP \left( \hty_1 =  t_N  \right) \qquad \text{ with } \eta'_{\gep} \stackrel{\gep\to 0}{\to} 0 \, .
\end{equation}

\smallskip
\noindent
{\bf Second step.}
It remains to control 
\begin{align}
\hP\Big( \cM_{1,\kappa_N} &\leq (1-\gep) t_N, \htx_{\kappa_N} = N , \hty_{\kappa_N} = \gamma_c N +t_N \Big) \notag\\
&= \hP\Big( \cM_{1,\kappa_N} \in \left( \gep t_N, (1-\gep) t_N \right), \htx_{\kappa_N} = N , \hty_{\kappa_N} = \gamma_c N +t_N \Big) 
\label{eq:part2.1}\\
& \quad + \hP\Big( \cM_{1,\kappa_N}\le \gep  t_N , \htx_{\kappa_N} = N , \hty_{\kappa_N} = \gamma_c N +t_N \Big)\, . 
\label{eq:part2.2}
\end{align}
The first term in the right-hand side, that is \eqref{eq:part2.1}, is smaller than
\begin{equation}
\begin{split}
&\sum_{k=1}^{N}  k \sum_{z=  \gep t_N  }^{  (1-\gep) t_N  } \sum_{x\in\bbN}  \hP \Big( \hty_1=z, \htx =x \Big) \bP \Big( \hty_{k-1} = \gamma_c N+ t_N -z , \htx_{k-1} = N-x \Big)\\
&\leq c \gep^{(2+\ga)} N   \bP \Big(\hty_1 =  t_N  \Big) \sum_{x\in\bbN}  p_h(x)  \sum_{k=1}^N  \bP \Big( \hty_{k-1} \geq  \gamma_c N + \gep t_N , \htx_{k-1} = N-x \Big)
\end{split}
\end{equation}
where we used \eqref{eq:uniformbound} uniformly for $z\geq \gep t_N$ and then summed over $z$ to get the first inequality. Then, we split the last sum into two parts. For $k\leq k_N^{(\gep)}:= N/\hat \mu_1 + \gep^2 t_N$, we have
\begin{equation}
\sum_{k=1}^{k_N^{(\gep)}}  \bP \left( \hty_{k-1} \geq  \gamma_c N + \gep t_N , \htx_{k-1} = N-x\right) \leq \hP \Big( \hty_{k_N^{(\gep)}}  \geq \gamma_c N +\gep t_N \Big).
\end{equation}
Then, provided that $\gep$ has been fixed small enough so that $\gamma_c N +\gep t_N \geq \hat \mu_2 k_N^{(\gep)} +\tfrac12 \gep t_N$, and since $t_N/a_N^{(2)}\to +\infty$ (and $a_{k_N^{(\gep)}}^{(2)} \leq a_N^{(2)}$), we have
\begin{equation}
\limsup_{N\to\infty}\hP \Big( \hty_{k_N^{(\gep)}}  \geq \gamma_c N +\gep t_N \Big) =0 \, . 
\end{equation}
On the other hand, for $k\geq k_N^{(\gep)}$, we have
\begin{equation}  
\sum_{k=k_N^{(\gep)} +1 }^{N} \hP \Big( \htx_{k-1} = N-x \Big)  \leq  \hP \Big( \htx_{k_N^{(\gep)}} \leq N -x \Big)  \leq  \hP \Big( \htx_{k_N^{(\gep)}} \leq \hat \mu_1 k_N^{(\gep)} - \hat \mu_1 \gep^2 t_N \Big),
\end{equation}
and since $t_N/\sqrt{N} \to +\infty$, also this terms goes to $0$ as $N\to\infty$.
In the end, we get that the term \eqref{eq:part2.1} is negligible compared to $N \hP \big( \hty_1 =  t_N  \big)$.

\smallskip
Then, it remains to bound \eqref{eq:part2.2}, and a first observation is that we can restrict it to having $\kappa_N  \leq k_N^+:= N/ \hat \mu_1 + t_N/ (4\hat \mu_2 ) $.
Indeed, we have that
\begin{equation}
\label{eq:forC0_1}
\begin{split}
 \hP \Big( \cM_{1,\kappa_N} \leq &\gep t_N , \htx_{\kappa_N}=N, \, \hty_{\kappa_N} = \gamma_c N +t_N , \kappa_N \geq k_N^+ \Big) \\
 &\le \hP \Big( \htx_{\kappa_N}=N, \kappa_N \geq k_N^+ \Big) \leq \hP \Big( \htx_{k_N^+ } \leq N \Big)\\
 & \le \exp\Big( - c\, (N- \hat \mu_1 k_N^+)^2 /k_N^+\Big) \le \exp\Big( - c'\, (t_N)^2 /N\Big),
\end{split}
\end{equation}
which decays faster than $N \hP \big( \hty_1=  t_N \big)$ because of assumption \eqref{hyp:bigjump}, provided that $C_0$ had been chosen large enough.

\smallskip
It remains to control
 \begin{align}
\sum_{k= 1 }^{k_N^+}   \hP \Big( \cM_{1,k} \leq \gep t_N , \htx_{k}=N, \, \hty_{k} = \gamma_c N +t_N \Big)   \,.
 \end{align}
We write that each term in the sum is
 \begin{equation}
 \sum_{j=  \log_2(1/\gep)}^{q_n}  \hP\left( \cM_{1,k} \in (2^{-(j+1)} t_N , 2^{-j} t_N ] , \htx_{k}=N, \, \hty_{k} = \gamma_c N +t_N \right)   \,,
 \end{equation}
 where 
$q_N$ is the  smallest integer such that $2^{-(q_N+1)} t_N <1$, so $q_N=O(\log_2N)$.
%,  and  we assume for conciseness that $ \log_2(1/\gep)\in \bbN$.
Then,  using \eqref{eq:uniformbound}, each term in the sum (i.e. for every $k$ and $j$) is bounded by a constant (not depending on $j$ and $k$) times
\begin{multline}
k \sum_{z =2^{-(j+1)} t_N }^{2^{-j} t_N } \sum_{x\in \bbN} \frac{L( 2^{-j} t_N)}{(2^{-j} t_N)^{(2+\ga)}} p_h(x)\times
\\
 \hP\left(\htx_{k-1}=N-x, \hty_{k-1} = \gamma_c N +t_N -z , \cM_{1,k-1}\leq 2^{-j} t_N \right) \\
\leq \frac{N  L(t_N) }
{t_N^{2+\ga} } \sum_{x\in \bbN} p_h(x) 2^{j(3+\ga)} \hP \bigg(\htx_{k-1}=N-x, \hty_{k-1} \geq \gamma_c N +\frac{t_N}2 , \cM_{1,k-1}\leq 2^{-j} t_N \bigg)\,,
\end{multline}
where we used that provided that $t_N$ is large enough, $L(2^{-j} t_N) \leq 2^{j} L(t_N)$ (this is a direct consequence of Potter's bound for slowly varying functions \cite[Th.~1.5.6]{cf:RegVar}) and summed over $z$.
Recovering the sum over $k$ and $j$, we therefore need to show that
\begin{align}
\label{eq:decompj}
\sum_{k=1}^{k_N^+} \sum_{j=\log_2(1/\gep)}^{q_n} \sum_{x\in\bbN} p_h(x) 2^{j(3+\ga)} \hP\Big(\htx_{k-1}=N-x, \hty_{k-1} \geq \gamma_c N +t_N/2 , \cM_{1,k-1}\leq 2^{-j} t_N \Big)
\end{align}
is small for $N$ large.

Then, for every $j$, we define $\{{\bar \tau}_k\}_{k= 0, 1, \ldots}$ (with distribution $\bar \bP^{(j)}$, carrying the  dependence on $j$) as an i.i.d.\ sum of $k$ variables with distribution $(\htx_1, \hty_1 \ind_{\{\hty_1 \leq 2^{-j }t_N\}})$: we therefore obtain that for $k\leq k_N^+$
\begin{equation}
\begin{split}
\hP &\Big(\htx_{k-1}=N-x, \hty_{k-1} \geq  \gamma_c N +t_N/2 ,  \cM_{1,k-1}\leq 2^{-j} t_N \Big) \\
&\leq \bar \bP^{(j)}\Big(\bar\tau^{(1)}_{k-1}=N-x, \bar\tau^{(2)}_{k-1} \geq \gamma_c N +t_N/2 \Big)
\leq \bar \bP^{(j)}\Big(\bar\tau^{(1)}_{k-1}=N-x, \bar\tau^{(2)}_{k_N^+} \geq \gamma_c N +t_N/2 \Big)\, .
\end{split}
\end{equation}
Using this inequality and summing it over $k$ in \eqref{eq:decompj}, (and then using that $\sum_x p_h(x) =1$), we obtain that \eqref{eq:decompj} is smaller than
\begin{align}
\label{eq:sumj}
\sum_{j=\log_2(1/\gep)}^{q_n} 2^{j(3+\ga)} \bar \bP^{(j)}\Big(\bar\tau^{(2)}_{k_N^+} \geq \hat \mu_2 k_N^+ +t_N/4 \Big) \, ,
\end{align}
where we used that $\gamma_c N \ge  \hat \mu_2 k_N^+ -\tfrac14 t_N$.
Then, we may use a Fuk-Nagaev inequality, see for example in~\cite{cf:Nag79}, to control the last probability -- we regroup the inequalities we need under  the following Lemma.

\begin{lemma}
\label{th:FukNag}
Let  $\{X_i\}_{i=1,2, \ldots}$ be a sequence of i.i.d.\ non negative  r.v.\ with $\bP(X_1>x) \sim \varphi(x) x^{-\rho}$ with $\rho>1$ and $\varphi(\cdot)$  a slowly varying function. Denote $\mu:=\bE[X]$ and $\sigma(y) = \bE[X_1^2 \ind_{\{X_1\leq y\}}]$. 
We have that there exists a constant $c>0$ such that for any $y\leq x$
\begin{equation}
\bP \bigg( \sum_{i=1}^n X_i \ind_{\{X_i\leq y\}} - \mu n \geq x \bigg) \leq 
\Big( c \frac{ n y^{1-\rho} \varphi(y) }{ x }  \Big)^{ c x/y }   +
 e^{-c x^2/n \sigma(y)}\ind_{\{\rho\geq 2\}}.
\end{equation}
\end{lemma}

\noindent
Lemma~\ref{th:FukNag} is taken from 
 \cite[Theorem 1.2]{cf:Nag79}  ($\rho\in(1,2)$) and      \cite[Corollary 1.6]{cf:Nag79} ($\rho\geq 2$).

Applying this lemma to $X_1 = \hty_1$ (i.e.
 $\rho=1+\alpha$ and $\varphi(\cdot)$ a constant times $L(\cdot)$), and $x= \gep t_N$, $y= 2^{-j} t_N$, we get that (using that $\sigma(y)\le \sigma(x)$ and $k_N^+ \leq N$ for the term $\ga\geq 1$)

\begin{equation}
\label{eq:termj}
\bar \bP^{(j)}\Big(\bar\tau^{(2)}_{k_N^+} \geq \mu_2 k_N^+ +t_N/4 \Big) \leq \left( \frac{c}{\gep} 2^{-j \ga} N (t_N)^{-(1+\ga)} L\left( 2^{-j} t_N\right) \right)^{c \gep 2^j} + e^{ - c'' \frac{t_N^2}{N \sigma(t_N)}  } \ind_{\{\ga\geq 1\}} \, .
\end{equation}
For the first term, we use that for $N$ large enough, $2^{-j \ga} L\left( 2^{-j} t_N\right) \leq L(t_N)$, so that it is bounded by
\begin{equation}
\left( \frac{c}{\gep}  N (t_N)^{-(1+\ga)} L\left( t_N\right) \right)^{c \gep 2^j}  = \Big( \frac{c}{\gep} N \hP\left( \hty_1 >t_N\right) \Big)^{c\gep 2^j}\leq  e^{ - 2^j }
 \end{equation}
where the second inequality holds for $N$ large enough (how large depends on $\gep$), since $N \hP \big( \hty_1 >t_N \big)\to 0$ as $N\to\infty$, simply because $t_N/m_N^{(2)} \geq c t_N/a_N^{(2)} \to +\infty$, recall \eqref{def:m_n}.

In the end, summing \eqref{eq:termj} for $j \in [\log_2(1/\gep), q_n]$, we get that \eqref{eq:sumj} is bounded by
\begin{equation}
\label{eq:forC0_2}
\sum_{j = \log_2(1/\gep) }^{q_n}  2^{j(3+\ga)} \Big( e^{-2^j} +  e^{- c'' t_N^2/ N \sigma(t_N)}  \ind_{\{\ga\geq 1\}} \Big)
 \leq  e^{- c/\gep} +  (q_n) N^{3+\ga} e^{- c'' t_N^2/ N \sigma(t_N)}  \ind_{\{\ga\geq 1\}}\, ,
\end{equation}
and the second term is small when $N\to +\infty$ thanks to assumption \eqref{hyp:bigjump}, provided that the constant $C_0$ has been fixed large enough in the case $\ga\geq 1$.
\qed

\section{The free partition function: proof of Theorem~\ref{th:aimZf}}
\label{sec:free}

%Suppose that $M = \gamma_N N $ with $\gamma_N \sim \gamma $ and $h>0$. In view of \eqref{eq:Zchat}, the sharp estimates of $Z^c_{N,M,h}$ are then reduced to use the estimates of the renewal function $\bP \left( (N,M) \in \hat{\tau}_h \right)$.  Since we work with a discret renewal process, all values in the summation in \eqref{eq:zf} have to be thought as the integer part, but we omit this in the notation to keep things simpler.
 
%To keep things simpler in this section, with some abuse of notation, we will systematically omit the integer part in the formulas. 
We will first prove \eqref{eq:aimZf1} for the case $\sum_j K_f(j) < \infty$. Many estimates are in common with the case 
 $\sum_j K_f(j) =\infty$ that we treat right after, and we will stress along the proof when the estimates are dependent or not on the fact that $\sum_j K_f(j) < \infty$.
 Also, the proof of the lower bounds \eqref{eq:aimZf2.1}
 and \eqref{eq:aimZf2.2} are  contained in the proof of \eqref{eq:aimZf1} as we explain along the way.

%We start by introducing the notations:
%\begin{equation}
%\begin{split}
%V:= "\text{length of the free part in the }y-\text{direction}",\\
%W:= "\text{length of the free part in the }x-\text{direction}",
%\end{split}
%\end{equation}
%so that we interpret in \eqref{eq:zf} the terms $K_f(i) = \bP(W=i)$ $K_f(j) =\bP(V=j)$.

\smallskip

\noindent
\emph{Proof of Theorem~\ref{th:aimZf}.}
As announced, we start with the proof  of \eqref{eq:aimZf1} and assume $\sum_j K_f(j) < \infty$.
Let us fix $\eta>0$, and $\gep>0$ small, how small depends on $\eta$ as will be stressed in the proof.

We decompose the free partition function into several parts: 
\begin{equation}
\label{eq:several}
\begin{split}
Z_{N,M,h}^{f} &= \mathrm{I} + \mathrm{II} + \mathrm{III} + \mathrm{IV} + \mathrm{V} \\
\text{with }\qquad
\mathrm{ I}&= Z_{N,M,h}^{f} \big( V_1^{(N)} \geq  1/\gep \big)\\
\mathrm{II}& = Z_{N,M,h}^f \big(V_1^{(N)}< 1/\gep ; V_2^{(N)} \leq  1/\gep \big) \\
\mathrm{III} & = Z_{N,M,h}^f  \big( V_1^{(N)} < 1/\gep ; V_2^{(N)} \in ( 1/\gep, t_N - \tfrac1\gep a_N^{(2)}) \big ) \\
\mathrm{IV} & = Z_{N,M,h}^f \big(V_1^{(N)}< 1/\gep ; V_2^{(N)} \in [t_N-\tfrac1\gep a_N^{(2)},t_N+\tfrac1\gep a_N^{(2)}] \big) \\
\mathrm{V} & = Z_{N,M,h}^f \big (V_1^{(N)}< 1/\gep  ; V_2^{(N)} >t_N+\tfrac1\gep a_N^{(2)} \big). 
\end{split}
\end{equation}
The main contribution comes from the terms $\mathrm{II}$ and $\mathrm{IV}$. We first estimate these terms, before showing 
that all the other ones are negligible compared to $\max(\mathrm{II},\mathrm{IV})$.

\subsection{Main terms, and proof of \eqref{eq:aimZf2.1} and \eqref{eq:aimZf2.2}}
\label{sec:free1}

\subsubsection{Analysis of II and proof of \eqref{eq:aimZf2.1}} 
This term can be written as
\begin{equation}
\mathrm{II}:= \sum_{i< 1/\gep } \sum_{j\leq 1/\gep} K_f(i) K_f(j) e^{(N-i) \tn_h} \bP( (N-i,M-j) \in\hat \tau) ,
\label{eq:useaimlb}
\end{equation}
and it is just a matter of estimating $\bP( (N-i,M-j) \in\hat \tau)$ uniformly for $0\leq i,j \leq 1/\gep$. 
We have from Theorem~\ref{th:aimZc}, uniformly for $i,j \leq 1/\gep$,
\begin{equation}
\bP( (N-i,M-j) \in\hat \tau) \stackrel{N \to \infty}\sim \frac{N-i}{\hat \mu_1^2 } \bP(\hat\tau^{(2)}_1 = M-j - \gamma_c(N-j)) \sim \frac{N}{\hat \mu_1^2} \bP(\hat\tau^{(2)}_1 = t_N)\, .
\end{equation}

Hence,
\begin{align}
\mathrm{II} & \sim  \frac{N}{\hat \mu_1^2} \bP(\hat\tau^{(2)}_1 = t_N) e^{N\tn_h }\sum_{i \leq 1/\gep} e^{-i \tn_h} K_f(i)\sum_{j\leq 1/\gep} K_f(j) \, .
\label{asympII}
\end{align}

Hence, if $\overline{ K}= \sum_{j\leq 1/\gep} K_f(j)  <+\infty$, we get that provided that $\gep$ has been fixed small enough (depending on $\eta$), for all $N$ sufficiently large,
\begin{equation}
\begin{split}
\mathrm{II} \ge (1-\eta) \overline{K} \, \frac{N}{\hat \mu_1^2}  e^{N\tn_h } \Big(\sum_{i \ge 0} e^{-i \tn_h} K_f(i) \Big) \bP(\hat\tau^{(2)}_1 = t_N) \, , \\
 \mathrm{II} \le (1+\eta) \overline{K} \, \frac{N}{\hat \mu_1^2}  e^{N\tn_h }  \Big(\sum_{i \ge 0} e^{-i \tn_h} K_f(i) \Big) \bP(\hat\tau^{(2)}_1 = t_N) \, .
  \end{split}
\end{equation}

The lower bound \eqref{eq:aimZf2.1} is obtained simply by using the estimate \eqref{eq:aimlb} instead of \eqref{eq:hatp} in \eqref{eq:useaimlb}: the straightforward details are left to the reader.  
%We refrained from doing so to avoid lengthy notations.

\subsubsection{Analysis of IV and proof of \eqref{eq:aimZf2.2}} 
It can be written as
\begin{equation}
\mathrm{IV} = \sum_{i\leq 1/\gep} \sum_{j=t_N -\tfrac1\gep a_N^{(2)}}^{t_N+\tfrac1\gep a_N^{(2)}} K_f(i) K_f(j) e^{(N-i) \tn_h} \bP( (N-i,M-j) \in\hat \tau).
\end{equation}
We have that uniformly for $j\in [t_N-\tfrac1\gep a_N^{(2)}, t_N+\tfrac1\gep a_N^{(2)}]$, $K_f(j) \sim K_f(t_N)$. We can therefore focus on estimating, uniformly for $i\leq 1/\gep$
\begin{align*}
\sum_{j=t_N -\tfrac1\gep a_N^{(2)}}^{t_N+\tfrac1\gep a_N^{(2)}}\bP\big( &(N-i,M-j) \in\hat \tau \big) \\
&= \bP\Big(\text{for some }k,\  \hat\tau^{(1)}_k = N-i , \hat\tau^{(2)}_k \in [\gamma_c N -\tfrac1\gep a_N^{(2)}, \gamma_c N +\tfrac1\gep a_N^{(2)}]  \Big)\, ,
\end{align*}
and now we prove that this term is close to $1/{\hat \mu_1}$.
In fact 
we have
\begin{multline*}
\bP\Big(\text{for some } k, \ \hat\tau^{(1)}_k = N-i ,  \hat\tau^{(2)}_k \in [\gamma_c N-\tfrac1\gep a_N^{(2)}, \gamma_c N+\tfrac1\gep a_N^{(2)}]\Big) \\
= \bP( N-i \in \hat \tau^{(1)}) - \bP\Big(\text{for some } k, \ \hat\tau^{(1)}_k = N-i , \hat\tau^{(2)}_k \notin [\gamma_c N-\tfrac1\gep a_N^{(2)}, \gamma_c N+\tfrac1\gep a_N^{(2)}] \Big),
\end{multline*}
and we now show that, provided that $\gep$ had been fixed small enough, uniformly for $i\leq 1/\gep$ and $N$ large enough:
\begin{align}
\bP\Big( \hat\tau^{(1)}_k = N-i , \hat\tau^{(2)}_k < \gamma_c N-\tfrac1\gep a_N^{(2)} \text{ for some } k \Big) \le \eta\, , \label{Ytoosmall}\\
\bP\Big( \hat\tau^{(1)}_k = N-i , \hat\tau^{(2)}_k >\gamma_c N+\tfrac1\gep a_N^{(2)} \text{ for some } k \Big) \le \eta \, . \label{Ytoolarge}
\end{align}
This will be enough, since by the Renewal Theorem we have that $\bP( N-i \in \hat\tau^{(1)}) \to \hat \mu_1^{-1}$ uniformly for $i\leq 1/\gep$.

\smallskip
 To treat \eqref{Ytoosmall}, define $k_N:= \frac{1}{\hat \mu_1 } N -  \tfrac{1}{2 \hat \mu_2 \gep} a_N^{(2)} $: we have uniformly for $i\le 1/\gep$
\begin{equation}
\begin{split}
 \bP \Big(\text{for } &\text{some } k, \ \hat\tau^{(1)}_k = N-i , \hat\tau^{(2)}_k < \gamma_c N-\tfrac1\gep a_N^{(2)} \Big)\\
& = \bP\Big( \text{for some } k \leq k_N,\ \hat\tau^{(1)}_k = N-i , \hat\tau^{(2)}_k < \gamma_c N-\tfrac1\gep a_N^{(2)}   \Big)  \\
& \qquad \quad+ \bP\Big(\text{for some } k \geq  k_N ,\ \hat\tau^{(1)}_k = N-i , \hat\tau^{(2)}_k < \gamma_c N-\tfrac1\gep a_N^{(2)}    \Big)\\
 &\leq \bP\left( \hat\tau^{(1)}_{k_N} \geq N - 1/\gep \right) + \bP\left( \hat\tau^{(2)}_{k_N} < \gamma_c N - \tfrac1\gep a_N^{(2)} \right).
\end{split}
\end{equation}
Now, it is easy to see that the two terms in the last line are small:
we indeed have that for arbitrary $\eta'>0$, one can choose $\gep$ small enough so that for all $N$ large enough,
\begin{equation}
\label{eq:gotozero}
\begin{split}
&\bP\left( \hat\tau^{(1)}_{k_N} \geq N - 1/\gep \right)  =\bP \left(\hat\tau^{(1)}_{k_N} \geq \hat \mu_1 k_N +\frac{\hat\mu_1}{2 \hat \mu_2 \gep} a_N^{(2)}  - 1/\gep \right) \le \eta' ,\\
\text{and }\quad & \bP\left( \hat\tau^{(2)}_{k_N} < \gamma_c N - \frac{1}{\gep} a_N^{(2)} \right)  = \bP \left(\hat\tau^{(2)}_{k_N} < \hat\mu_2 k_N -  \frac{1}{2\gep} a_N^{(2)} \right) \le \eta' \, .
\end{split}
\end{equation}
For the first line, we used that $\frac{\hat\mu_1}{2 \hat \mu_2 \gep} a_N^{(2)}  - 1/\gep  \ge \gep^{-1/2} \sqrt{N} \ge \gep^{-{1/2}} \sqrt{k_N}$ provided that $N$ is large enough (and $\gep$ small), and then simply Chebichev's inequality. For the second line, we used that $ \gamma_c = \hat \mu_2/\hat\mu_1 $ to get that $\gamma_c N  = \hat\mu_2 k_N + \tfrac{1}{2\gep} a_N^{(2)}$, and then the approximation of $ (a_{k_N}^{(2)} )^{-1}(\hat\tau_{k_N} - k_N \hat \mu_2)$ by an $\alpha_2$-stable distribution, as done in \eqref{eq:est-E_3}.

\smallskip
For \eqref{Ytoolarge}, we define $k'_N = \frac{1}{\hat \mu_1 } N +  \frac{1}{2\hat \mu_2 \gep} a_N^{(2)}$, and similarly to what is done above, we have
\begin{align*}
\bP\Big(\text{for some } k,\ \hat\tau^{(1)}_k = N-i ,& \hat\tau^{(2)}_k > \gamma_c N + \tfrac1\gep a_N^{(2)}  \Big)\\
&\leq \bP\left( \hat\tau^{(1)}_{k'_N} \leq  N - i \right) + \bP\left( \hat\tau^{(2)}_{k'_N} >  \gamma_c N  + \tfrac1\gep a_N^{(2)}  \right) ,
\end{align*}
and both terms are smaller than $\eta'$ provided that $\gep$ had been fixed small enough and $N$ is large, for the same reasons as in \eqref{eq:gotozero}.

In the end, we get that provided that $\gep$ had been fixed small enough, for all sufficiently large $N$
\begin{equation}
\begin{split}
\mathrm{IV} &\ge  (1-\eta) \frac{1}{\hat \mu_1} K_f(t_N)   e^{N\tn_h } \sum_{i\le 1/\gep} e^{-i \tn_h} K_f(i) \, ,\\
\mathrm{IV} & \le  (1+\eta) \frac{1}{\hat \mu_1} K_f(t_N)   e^{-N\tn_h }  \sum_{i \le 1/\gep} e^{-i \tn_h} K_f(i)\, .
\label{asympIV}
\end{split}
\end{equation}
Obviously, since the last sum converges, we can replace it with the infinite sum, and simply replace $\eta$ by $2\eta$ provided that $\gep$ is small enough. This competes the analysis of IV.

\smallskip

For what concerns  \eqref{eq:aimZf2.2} we simply need to show that 
\begin{align*}
\mathrm{IVb}&:= Z_{N,M,h}^f \Big( V_1^{(N)}\le \frac1\gep, V_2^{(N)} \in \big[t_N -\tfrac1\gep a_N^{(2)} , t_N+ \tfrac1\gep a_N^{(2)}\big], \mathcal{M}_{1,\kappa_N} > \tfrac1\gep m_N^{(2)} \Big)\\
&= \sum_{i\leq 1/\gep} \sum_{j=t_N -\tfrac1\gep a_N^{(2)}}^{t_N+\tfrac1\gep a_N^{(2)}} K_f(i) K_f(j) e^{(N-i) \tn_h} \bP\Big( (N-i,M-j) \in\hat \tau, \mathcal{M}_{1,\kappa_N} >\tfrac1\gep m_N^{(2)}  \Big).
\end{align*}
is negligible compared to \eqref{asympIV}. But again, uniformly for the range of $j$ considered, we have $K_f(j) \le 2 K_f(t_N)$ (provided
that $N$ is large enough). Then, dropping the event $N-i\in \htx$, and summing over $j$, we get that
\begin{align*}
\mathrm{IVb} \le 2 K_f(t_N) e^{N \tn_h} \Big(\sum_{i\le 1/\gep} K_f(i) e^{-i \tn_h} \Big) \bP\Big(  \mathcal{M}_{1,\kappa_N} >\tfrac1\gep m_N^{(2)} \Big)\, .
\end{align*}
Then, using that $\kappa_N \le N$, we get that
\begin{align*}
\bP\Big(\mathcal{M}_{1,\kappa_N} >\tfrac1\gep m_N^{(2)}  \Big) \le \bP\Big( \max_{1\le k\le N} (\hty_k -\hty_{k-1})>\frac1\gep m_N^{(2)}   \Big) \le N \bP\big( \hty_1 > \tfrac1\gep m_N^{(2)}\big), 
\end{align*}
which can be made arbitrarily small by choosing $\gep$ small (uniformly in $N$), thanks to the definition \eqref{def:m_n} of $m_N^{(2)}$.
Hence $\mathrm{IVb}$ is negligible compared to $\mathrm{IV}$.
We also stress here that to estimate $\mathrm{IV}$ -- in particular to obtain \eqref{asympIV}~--, we did not make use of the assumption $\sum K_f(i) <+\infty$.

\subsection{Remaining terms}
\label{sec:free2}

It remains to estimate the terms $\mathrm{I}$, $\mathrm{III}$ and $\mathrm{V}$ in \eqref{eq:several}, and show that they are negligible compared to \eqref{asympII} or \eqref{asympIV}.
We start by parts $\mathrm{III}$ and $\mathrm{V}$.

\subsubsection{Analysis of III} 
Assume that $N$ is large enough, so that  $\tfrac1\gep a_N^{(2)} \leq \tfrac12 t_N$ we write
\begin{align}
\label{eq:III-1}
\mathrm{III}\leq & Z_{N,M,h}^f \Big( V_1^{(N)}< 1/\gep ; V_2^{(N)} \in ( 1/\gep,  \tfrac12 t_N)  \Big)  \quad (\text{denoted }\mathrm{IIIa}) \notag\\
&+ Z_{N,M,h}^f \Big( V_1^{(N)}< 1/\gep ; V_2^{(N)} \in ( \tfrac12 t_N, t_N - \tfrac1\gep a_N^{(2)})  \Big)  \quad (\text{denoted }\mathrm{IIIb})
\end{align}

The first term is
\begin{equation}
\begin{split}
\mathrm{IIIa} = \sum_{i <1/\gep } \sum_{j =1/\gep}^{t_N/2}   K_f(i) K_f(j) e^{(N-i) \tn_h} \bP \big( (N-i, M -j) \in\hat \tau \big) .
\end{split}
\end{equation}
Now we can bound, uniformly for $i < 1/\gep$ and $j\leq t_N/2$ (so that $(M-j) - \gamma_c(N-i) \geq t_N/4$ for $N$ sufficiently large, and we can apply Theorem \ref{th:aimZc})
\begin{equation}
\bP\big( (N-i, M -j) \in\hat \tau \big) \leq c N  \sup_{m\geq t_N/4 } \bP(\hat\tau^{(2)}_1 = m) \leq c' N  \bP(\hat\tau^{(2)}_1=t_N) .
\label{eq:boundproba}
\end{equation}
Hence we get
\begin{equation}
\label{eq:boundproba-2}
\mathrm{IIIa} \le
c' N \bP\left(\hat\tau^{(2)}_1 =t_N\right) e^{N \tn_h}\sum_{i} e^{-i \tn_h} K_f(i) \sum_{j =1/\gep}^{t_N/2} K_f(j),
\end{equation}
and in the case when $\sum K_f(j) <+\infty$, the last sum can be made arbitrarily small by choosing $\gep$ small.
Hence, recalling \eqref{asympII}, we get that $\mathrm{IIIa}\le \eta \times \mathrm{II}$ for all $N$ sufficiently large, provided that $\gep$ is small enough.

For the term $\mathrm{IIIb}$, we use that $K_f(j) \leq c K_f(t_N)$ uniformly for $j\geq t_N/2$, to get that
\begin{multline}
\mathrm{IIIb} \le
\sum_{i\leq 1/\gep} \sum_{j =t_N/2}^{t_N - \tfrac1\gep a_N^{(2)}}   K_f(i) c K_f(t_N) e^{(N-i) \tn_h} \bP\left( (N-i, M -j) \in\hat \tau\right)
\\
\le  c\, K_f(t_N) e^{N \tn_h}\sum_{i\leq 1/\gep} e^{-i \tn_h} K_f(i) \bP\big( \text{for some }k,\ \hat\tau^{(1)}_k = N-i , \hat\tau^{(2)}_k > \gamma_c N + \tfrac1\gep a_N^{(2)} \big).
\end{multline}
Since we have seen in \eqref{Ytoolarge} that the last probability is smaller than some arbitrary $\eta'$ for all $N$ large enough (provided that $\gep>0$ is small enough), uniformly for all $i\le 1/\gep$ we have that
$\mathrm{IIIb}  \le \eta\times  \mathrm{IV}$ (recall \eqref{asympIV}).
We stress that, here again, we do not make use of the assumption $\sum K_f(i)<+\infty$.

\smallskip
In the end, we obtain that $\mathrm{III} \le \eta \times (\mathrm{II}+  \mathrm{IV})$ (provided that $N$ is large enough).

\subsubsection{Analysis of V}
We proceed analogously as above: using that $K_f(j) \leq c K_f(t_N)$ uniformly for $j\geq t_N$, we get that
\begin{equation}
\notag
\begin{split}
\mathrm{V} & \leq \sum_{i\leq 1/\gep} \sum_{j \geq t_N + \tfrac1\gep a_N^{(2)} }   K_f(i) c K_f(t_N) e^{(N-i) \tn_h} \bP\big( (N-i, M -j) \in\hat \tau \big) \\
& \leq c K_f(t_N) e^{N \tn_h}\sum_{i\leq 1/\gep} e^{-i \tn_h} K_f(i) \bP\big( \text{for some } k,\ \hat\tau^{(1)}_k = N-i , \hat\tau^{(2)}_k < \gamma_c N - \tfrac1\gep a_N^{(2)} \big).
\end{split}
\end{equation}
Now we again recall \eqref{Ytoosmall}, which tells that the last probability is smaller than some arbitrary $\eta'$ for all $N$ large enough (provided that $\gep>0$ is small enough, uniformly for all $i\le 1/\gep$). In the end, in view of \eqref{asympIV}, we get that $\mathrm{V} \le \eta \times \mathrm{IV}$, and here again we did not make use of the assumption $\sum K_f(i)<+\infty$.

\subsubsection{Analysis of I} 
We separate it into two parts: $V_1^{(N)} \geq (\log N)^2$, and  $V_1^{(N)} \in (1/\gep , (\log N)^2)$. 
We have
\begin{align}
\label{eq:I1b}
Z_{N,M,h}^{f} \Big(& V_1^{(N)}\geq (\log N)^2 \Big)   = \sum_{i \geq (\log N)^2} \sum_{j = 0}^M K_f(i) K_f(j) e^{(N-i) \tn_h} \bP\big( (N-i, M-j) \in \hat\tau \big) \notag\\
&\leq  \Big(\sum_{j = 0}^M K_f(j) \Big) e^{N \tn_h} \sum_{i = (\log N)^2}^N e^{- i \tn_h} K_f(i)  \leq  c N^{c+2} e^{N \tn_h} e^{- (\log N)^2 \tn_h},
\end{align}
where we first simply bounded the probability by $1$, and also that there is some constant $c<0$ such that $K_f(i)\le c N^{c} $  for $i\le N,M$.
Clearly,  in view of \eqref{asympIV} (or \eqref{asympII}), we get that
$Z_{N,M,h}^{f}(V_1^{(N)} \geq (\log N)^2 )  = o(\mathrm{IV})$, since $1/K_f(t_N)$ and $1/\bP(\hat\tau^{(2)}_1 = t_N)$ are  $O(N^{c'})$ for some $c'>0$. Again, we did not use that $\sum K_f(i)<+\infty$, even if it would have simplified the upper bound.

\smallskip
We now turn to the case when $V_1^{(N)} \le (\log N)^2$.
We write
\begin{equation}
\begin{split}
Z_{N,M,h}^{f}  \Big(V_1^{(N)} \in [1/\gep , (\log N)^2) \Big) \, & =\, 
Z_{N,M,h}^{f}\Big(V_1^{(N)} \in [1/\gep , (\log N)^2), V_2^{(N)} \leq t_N/2\Big)\, 
\\
 & + Z_{N,M,h}^{f} \Big(V_1^{(N)} \in [1/\gep , (\log N)^2), V_2^{(N)} > t_N/2 \Big)\, .
\end{split}
\end{equation}

For the first term, and using that $\bP\big( (N-i,M-j) \in\hat \tau \big) \leq c N \bP( \hat\tau^{(2)}_1 = t_N)$ uniformly for $i\leq (\log N)^2$ and $j\leq t_N/2$ (since then we have $M-j-N-i \ge t_N/4$ for $N$ large enough, similarly to \eqref{eq:boundproba}), we have
\begin{multline}
\label{eq:I1b-2}
Z_{N,M,h}^{f} \Big( V_1^{(N)} \in [1/\gep , (\log N)^2), V_2^{(N)} \leq t_N/2 \Big) \\
\leq \sum_{i\geq 1/\gep} \sum_{j=0}^{t_N/2} K_f(i) K_f(j)e^{(N-i) \tn_h} c N \bP(\hat\tau^{(2)}_1=t_N)  \\
\leq c \Big( \sum_{j=0}^{t_N/2} K_f(j) \Big) N \bP(\hat\tau^{(2)}_1 =t_N) e^{N \tn_h} \sum_{i\ge 1/\gep} K_f(i) e^{- i \tn_h}\, .
\end{multline}
When $\sum_j K_f(j) <+\infty$, then recalling \eqref{asympII}, this term is smaller than $\eta \times \mathrm{II}$ provided that $\gep$ is small enough.

For the second term, we use that $K_f(j)\leq c K_f(t_N)$ uniformly for $j\geq t_N/2$ to get 
\begin{multline}
\label{eq:I1b-3}
Z_{N,M,h}^{f} \Big( V_1^{(N)} \in (1/\gep , (\log N)^2), V_2^{(N)} > t_N/2 \Big) \\
\leq \sum_{i\geq 1/\gep} \sum_{j=t_N/2 +1}^{M} K_f(i) c K_f(t_N) e^{(N-i) \tn_h} \bP\big( (N-i,M-j)\in\hat \tau \big) \\
\leq c K_f(t_N) e^{N \tn_h} \sum_{i\ge 1/\gep} K_f(i) e^{- i \tn_h} \, ,
\end{multline}
where we used that the sum over $j$ of $\bP\big( (N-i,M-j)\in\hat \tau \big)$ is bounded by $1$. Then, the last sum can be made arbitrarily small by choosing $\gep$ small, so that in view of \eqref{asympIV}, this term can be bounded by $\eta \times \mathrm{IV}$ (and note that we did not use that $\sum_j K_f(j)< \infty$). 

\subsubsection{Conclusion in the case of $\sum_j K_f(j)< \infty$}
We have therefore proven that for any $\eta>0$, we can choose $\gep>0$ small such that, for all $N$ large enough (how large depend on $\gep$), 
\begin{equation}
\mathrm{II} + \mathrm{IV} \le Z_{N,M,h}^f  \le (1+2\eta) \mathrm{II} + (1+4\eta)\mathrm{IV}
\end{equation}
and the two terms behave asymptotically respectively as \eqref{asympII} and \eqref{asympIV}: this proves \eqref{eq:aimZf1} for $\sum_j K_f(j)< \infty$.

%
%Moreover,
%since $u_N$ and $v_N/a^{(2)}_N$ may be chosen to diverge arbitrarily slowly, it is elementary to check that the estimates we just performed
%hold also if we perform by fixing $u_N=1/\gep$ and $v_N/a^{(2)}_N=1/\gep$ and if we take the limit $\gep \searrow 0$ after
%$N \to \infty$. More explicitly we have:
%\begin{multline} 
%Z_{N,M,h}^{f} = \eta_1 (\gep, N) Z_{N,M,h}^f(V_1 < 1/\gep, V_2<1/\gep)\,  + \\
%\eta_1 (\gep, N) Z_{N,M,h}^f(V_1 < 1/\gep, V_2\in(t_N -a^{(2)}_N/ \gep  ,t_N+a^{(2)}_N/ \gep))\, , 
%\end{multline}
%with $\lim_{\gep \searrow 0} \limsup_{N \to infty}\vert \eta_i (\gep, N) -1\vert=0$ for $i=1,2$. This observation directly implies  
%\eqref{eq:aimZf2.1} and \eqref{eq:aimZf2.1} for $\sum_j K_f(j)< \infty$.

\subsubsection{The case of $\sum_j K_f(j)= \infty$, with $\overline{\ga} <1$}
This time we have to show that IV dominates. We go through the various terms, but as pointed out during the proof,
we have not used $\sum_j K_f(j)< \infty$ in estimating  IV,  so \eqref{asympIV} still holds.
We retain, for local use, that  IV behaves (and, in particular, is bounded from below by) a constant times  
$K_f(t_N)  \exp(N \tn_h)$.

The estimate \eqref{asympII} for II is still valid.
This term can be dealt directly without troubles, but it is more practical to observe that this time II is dominated by IIIa (for $N$ large, of course).  We can therefore focus on \eqref{eq:boundproba-2} which, up to a constant 
factor,  is bounded by
\begin{equation}
N \bP \Big( \hat\tau^{(2)}_1=t_N \Big) \exp(N \tn_h) \sum _{j \le t_N/2} K_f(j) \,.
\end{equation}
Therefore, in view of the behavior of IV that we have just recalled, this term is negligible  if
\begin{equation}
\label{eq:forIII}
N t_N  \bP \Big( \hat\tau^{(2)}_1=t_N \Big)  \, \ll\, 1\, ,
\end{equation} 
since  $\sum _{j \le t_N/2} K_f(j) \le cst. t_N K_f(t_N) $ if $\overline{\ga}<1$.

But the left-hand side is equivalent to 
$N L(t_N)/t_N^{1+\ga}$ and hence \eqref{eq:forIII} directly follows by recalling the definition \eqref{def:an}
of $a_N^{(2)}$ and that $t_N \gg a_N^{(2)}$.
This shows that both II and of IIIa are negligible compared to IV.

The estimates for IIIb and V, as already pointed out, are valid without assuming that $\sum_j K_f(j)< \infty$,
so we are left with controlling I.
Recall that we split the contribution of I into three parts: \eqref{eq:I1b}, \eqref{eq:I1b-2} and \eqref{eq:I1b-3}.
As noticed above, the fact that $\sum_j K_f(j)< \infty$ was not used in estimating \eqref{eq:I1b} and \eqref{eq:I1b-3}. Moreover \eqref{eq:I1b-2} we can be bounded like  \eqref{eq:boundproba-2} (in fact, it is much smaller), that was found above to be negligible compared to IV. We therefore conclude that I is also negligible compared to IV, and this completes the analysis of the
case $\sum_j K_f(j) =\infty$, and of the proof of Theorem \ref{th:aimZf}.
\qed
\medskip

\section*{Acknowledgements}
G.G. acknowledges the support of grant ANR-15-CE40-0020.

\begin{appendix}

\section{The case  $\overline{\ga} =1$ and $\sum_j K_f(j) =+\infty$}
\label{app:baralpha}

We treat this case in a concise way
because most of the technical work has already been done above.
To summarize -- recall the different contributions in \eqref{eq:several} -- the term IV is well estimated in \eqref{asympIV} and the terms I, IIIb and V were found to be negligible compared to it -- this was valid even when $\sum_j K_f(j) =+\infty$. When $\sum_{j} K_f(j) =+\infty$, then the term II is found to be negligible compared to IIIa, and we therefore focus on this last term.

We can again decompose IIIa into two contributions: 
\[\mathrm{IIIa} = Z_{N,M,h}^f \Big( V_1^{(N)}< 1/\gep ; V_2^{(N)} \in ( 1/\gep,  \gep t_N)  \Big)  + Z_{N,M,h}^f \Big( V_1^{(N)}< 1/\gep ; V_2^{(N)} \in (  \gep t_N , t_N/2)  \Big)  \, .\]
The second one, exactly in the same manner as for IIIb,  can be shown to be negligible compared to IV as $N\to\infty$.
Then, the first term is equal to
\[\mathrm{IIIa'}:=\sum_{i <1/\gep } \sum_{j =1/\gep}^{\gep t_N}   K_f(i) K_f(j) e^{(N-i) \tn_h} \bP \big( (N-i, M -j) \in\hat \tau \big) .\]
Then, thanks to Theorem \ref{th:aimZc}, for every $\eta>0$ we can choose $\gep>0$ small enough and $N_{\gep}$ large enough so that  uniformly for the range of $i$ and $j$ considered, and $N\ge N_{\gep}$
\[
\bP \big( (N-i, M -j) \in\hat \tau \big)
\begin{cases}
\ge (1-\eta) \frac{N}{\hat \mu_1^2} \bP(\hat\tau_1^{(2)} = t_N) \, , \\
\le (1+\eta) \frac{N}{\hat \mu_1^2} \bP(\hat\tau_1^{(2)} = t_N) \, ,
\end{cases}
\]
and we stress that the main contribution to this probability comes from a big loop event, of length larger than $(1-\gep) t_N$.
We therefore get that, for $N$ large enough, and denoting $\overline{K} (x) := \sum_{j=1}^x K_f(j)$ which is a slowly varying function, 
\begin{align*}
\mathrm{IIIa'}& \ge (1-\eta) \frac{N}{\hat \mu_1^2} \bP(\hat\tau_1^{(2)} = t_N)  e^{N\tn_h} \Big( \sum_{j =1/\gep}^{\gep t_N}K_f(j) \Big)  \sum_{i <1/\gep }  K_f(i)  e^{-i \tn_h}\\
&\ge (1-\eta)  \frac{N}{\hat \mu_1^2} \bP(\hat\tau_1^{(2)} = t_N) \overline{K}(t_N) e^{N \tn_h}  \Big(\sum_{i \ge 0}  K_f(i)  e^{-i \tn_h} \Big)\, ,
\end{align*}
and similarly for an upper bound with $1-\eta$ replaced by $1+\eta$.

We are actually able to narrow the condition $V_2^{(N)} \in (1/\gep, \gep t_N)$ in $\mathrm{IIIa'}$ to a smaller interval $( v_N , \gep_N t_N)$ without changing the estimates, provided that $v_N \to+\infty$ and $\gep_N \to 0$ slowly enough, precisely:
\begin{equation}
\label{cond:vNgepN}
\overline{K}(v_N) \ll \overline{ K}(\gep_N t_N) \, , \quad \quad \quad    \overline{ K}(\gep_N t_N) \stackrel{N\to\infty}{\sim} \overline{ K}(t_N) \, .
\end{equation}

We end up with the following result: recall the definition \eqref{eq:BLBS} of the Big Loop and Unbound strand, and when $\overline{\ga} =1$ with $\sum_j K_f(j)=+\infty$, define the new event $E_{mixed}^{(N)}$
\begin{equation}
E_{mixed}^{(N)} = \Big\{ \frac{\mathcal{M}_{1,\kappa_N}}{t_N} \in \big[  1-\gep_N , 1+\gep_N \big]  , \mathcal{M}_{2,\kappa_N} < m_N^+, V_1^{(N)} \leq u_N, V_2^{(N)} \in \big[v_N, \gep_N t_N \big] \Big\}
\end{equation}
where $v_N \gg 1$ and $\gep_N \ll 1$ are chosen as in \eqref{cond:vNgepN}. The event $E_{mixed}^{(N)}$ is therefore a set of trajectories with both a big loop (of order $t_N$), and a large unbound strand (of large order, but much smaller than $t_N$) -- to optimize the interval for the length of the unbound strand, one can take $v_N \to +\infty$ and $\gep_N \to 0$ as fast as possible, with the limitation given by \eqref{cond:vNgepN}.

\begin{theorem}
\label{th:baralpha}
Assume that $\ga >0$ and \eqref{hyp:bigjump1}, and if $\ga\geq 1$ assume additionally \eqref{hyp:bigjump}. We assume that $\overline{\ga}=1$ and that $\sum_j K_f(j) = +\infty$, and we denote $\overline{K}(x):= \sum_{j=1}^x K_f(j)$. 
Then, as $N\to\infty$, 
\begin{equation}
Z_{N,M, h}^{f} = (1+o(1)) Z_{N,M,h}^f \Big(E_{mixed}^{(N)}  \Big) + (1+o(1)) Z_{N,M,h}^f \Big(E_{US}^{(N)}  \Big) \, ,
\end{equation}
with
\begin{align}
e^{-N \tn_h}Z_{N,M,h}^f \Big(E_{mixed}^{(N)}  \Big) &\stackrel{N\to\infty}{\sim} \frac{N}{\hat \mu_1^2} \bP(\hat\tau_1^{(2)} = t_N)  \overline{ K}(t_N)  \Big(\sum_{i \ge 0}  K_f(i)  e^{-i \tn_h} \Big) \, , 
\label{asympEmix}\\
e^{-N \tn_h}Z_{N,M,h}^f \Big(E_{US}^{(N)}  \Big)  &\stackrel{N\to\infty}{\sim} \frac{1}{\hat \mu_1}\bigg(\sum_{i\geq 0} K_f(i) e^{- i \tn_h} \bigg)  K_f\left( t_N\right) \, .
\label{asympEUS}
\end{align}
\end{theorem}

Obviously, this theorem can easily be translated in term of path properties.
Indeed, since $\bP(\hat\tau_1^{(2)} = t_N) \sim cst.\, t_N^{-1} \bP(\hat\tau_1^{(2)} > t_N)$ and $K_f(t_N) = \overline{L}(t_N) t_N^{-1}$,  we have the asymptotic of the ratio 
\begin{equation}
\label{eq:tildeQ}
\tilde Q_N =\tilde Q_N(t_N):=\frac{Z_{N,M,h}^f \big(E_{mixed}^{(N)}  \big)}{Z_{N,M,h}^f \big(E_{US}^{(N)}  \big)}\stackrel{N\to\infty}{\sim} cst.\, N \bP(\hat\tau_1^{(2)} > t_N)  \frac{\overline{K} (t_N)}{ \overline{L}(t_N)},
\end{equation} 
with $\overline{K}(x) /\overline{L}(x) \to +\infty$ as a slowly varying function.
Therefore, we obtain
\begin{equation}
\label{eq:Palphabar}
\bP_{N,M,h}^f (E_{US}^{(N)}) \stackrel{N\to\infty}{\sim} \frac{1}{1+\tilde Q_N} \quad \text{and } \quad\bP_{N,M,h}^f (E_{mixed}^{(N)}) \stackrel{N\to\infty}{\sim}\frac{\tilde Q_N}{1+\tilde Q_N} \, .
\end{equation}

We stress that, when $\ga>1$, the ration $\tilde Q_N$ always goes to $0$ as $N\to \infty$: indeed, in that case $N \bP(\hat\tau_1^{(2)} > t_N) $ decays faster than any slowly varying function.
However, in the case $\ga\in (0,1]$, the ratio $\tilde Q_N$ diverges when $t_N \to +\infty$ slowly enough, showing that there is a regime under which the mixed trajectories described in the event $E_{mixed}^{(N)}$ occur, in the sense that $\bP_{N,M,h}^f (E_{mixed}^{(N)}) \to 1$ as $N\to\infty$.

\section{About the transition 
between  Cram\'er and\\ non-Cram\'er regimes}
\label{app}

In this Appendix, we discuss the condition \eqref{hyp:bigjump1}-\eqref{hyp:bigjump} ensuring that one lies in the big-jump regime described by Theorem~\ref{th:aimZc}. We focus on the constrained partition function -- or rather the probability $\bP((N,M)\in\hat \tau)$ -- to study the transition between the condensation phenomenon that we highlighted and the Cram\'er regime, but all the observations made here could also apply to the other results. Like in Section~\ref{sec:free} we omit integer parts, so $\gamma_c N$ stands for  the (upper or lower, as one wishes) integer part of $\gamma_c N$.

\subsection{Between Cram\'er and non-Cram\'er regimes I}
\label{app1}

If one sets $M=\gamma_c N$, or in other words if $t_N=0$, then \cite{cf:B} proves that 
\begin{equation}
\label{cramerborder}
\bP \big( (N,\gamma_c N) \in\hat \tau \big) \stackrel{N\to\infty}{\sim} \frac{c_0}{ a_N^{(2)}}\, ,
\end{equation}
 where the constant $c_0>0$ is explicit. The heuristics of this result can be easily understood: the typical number of renewal is $k_N = N/\hat\mu_1 +O( \sqrt{N})$ and, for each $k$ in that range, Doney's Local Limit Theorem \cite{cf:DonLLT} gives that $\bP(\hat\tau_k = (N,M))$
  is equivalent up to a multiplicative constant to $ ( a_N^{(2)} \sqrt{N})^{-1}$. Hence,  neither $\htx$ nor $\hty$ have to make an atypical deviation, and the term $(a_N^{(2)})^{-1}$ simply comes from a local limit theorem: there is no condensation phenomenon, 
  i.e. the typical trajectories contributing to the event $(N,\gamma_c N)\in\hat \tau$ do not exhibit a big jump. However, we are not 
  in the Cram\'er regime -- one component of the inter-arrivals does not have exponential tails, so there are jumps that are luch larger than $\log N$ -- and we can see this critical situation as a \emph{moderate} Cram\'er regime, because (moderate) deviations are carried by both components, like in the    Cram\'er regime the (large) deviations are carried by both components.
  
 A behavior like \eqref{cramerborder} also holds when $t_N/a_N \to t \in\bbR$: the constant $c_0$ is simply replaced by a constant $c_t$ depending on $t$. 
When $\ga \in (0,1)$, the fact that one lies in the big-jump regime (and Theorem \ref{th:aimZc} holds) as
soon as $t_N/a_{N}^{(2)}\to +\infty$ is optimal, in the sense that when $\sup_N t_N /a_N^{(2)} <+\infty$, then the typical trajectories do not exhibit a condensation phenomenon.

\subsection{Between Cram\'er and non-Cram\'er regimes II}
\label{app2}

When $\ga \ge 1$, the situation is more involved because the condition $t_N/a_N^{(2)}\to+\infty$ alone is not enough to  ensure that the model  is in the big-jump domain. 

We conjecture that when $\ga>1$, there is some $a_c =a_c(\ga)$ -- that we give explicitly below~-- such that the big-jump regime holds when $t_N > a \sqrt{N\log N}$ with $a>a_c$ (i.e.\ theorem \ref{th:aimZc} holds), and a moderate Cram\'er regime holds when $t_N < a \sqrt{N\log N}$ with $a<a_c$ (we give an explicit conjectured analogue of Theorem~\ref{th:aimZc}, see \eqref{eq:conj} below).
Finding the correct threshold when $\ga=1$ is even more involved and we prefer to leave it aside.
%: a weak statement would be that a sufficient to lie in the big-jump domain is that $t_{N}/a_N^{(2)} \geq a\sqrt{\log N}$ for some $a>0$ (note that the value $a_c$ found for $\ga>1$ goes to $0$ as $\ga\downarrow1$). 

So let us now focus on the case $\ga>1$, and develop some heuristic arguments to conjecture the asymptotic behavior of $\bP\big( (N,M)\in\hat \tau\big)$,
and the typical behavior of trajectories contributing to this event. We take $a_N^{(2)}=\sqrt{N}$, and we are considering the case $t_N/\sqrt{N} \to\infty$ (the case $t_N /\sqrt{N} \to t\in\bbR$ being given in Appendix \ref{app1}), with $t_N\le C_0 \sqrt{N\log N}$ (otherwise we already know we are in the big-jump domain).
Writing
\begin{equation}
\hP \big( (N,\gamma_c N +t_N) \in\hat\tau \big) = \sum_{k=1}^N \hP\Big( \htx_k= N, \hty_k= \gamma_c N +t_N\Big),
\end{equation}
then, the $k$'s bringing the main contribution to the sum are either  $k = N /\hat\mu_1 + O(\sqrt{N})$, in which case the deviation is entirely carried by $\hty$;  $k= N /\hat\mu_1 + t_N / \hat\mu_2+ O(\sqrt{N})$, in which case the cost is brought by $\htx$; or more generally  $k= N/\hat \mu_1 +  \theta t_N /\hat \mu_2  +O(\sqrt{N})$ with some $\theta\in\bbR$ (it is natural to expect $\theta\in[0,1]$, but $\theta\notin [0,1]$ should not be excluded), in which case the cost is shared \emph{jointly} by the two coordinates $\htx,\hty$.

Then, having a look at Nagaev's Theorem 1.9 in \cite{cf:Nag79} suggests that for any $k$, only two possible behavior can contribute to  $\hP ( \htx_k= N, \hty_k= \gamma_c N +t_N)$: having one large jump (in which case, and since $\hty$ has a heavier tail, the probability is maximal when $k =N/\hat\mu_1 + O(\sqrt{N})$ so that only $\hty$ has to make a large jump), or using a collective joint strategy with no big jump (i.e.\ a moderate Cram\'er regime).
The first possible behavior is therefore the big-jump strategy that we already identified, and we would therefore have that 
\begin{equation}
\label{eq:conj1}
\begin{split}
\hP\Big( (N,\gamma_cN & +t_N) \in\hat\tau\Big) = (1+o(1)) \frac{N}{\hat \mu_1^2} \hP\big( \hty_1=t_N\big)\\
& +  (1+o(1)) \sum_{k=1}^{N} \hP\big(\htx_k= N, \hty_k= \gamma_c N +t_N, \text{``with no big jump''}  \big)  \, ,
\end{split}
\end{equation}
where by ``no big jump'' we mean that all jumps are $O(m_N)$.

Using a local moderate deviation theorem for the probability when no big jump occurs (such a local moderate deviation theorem should hold because $t_N$ is not too large, $t_N \le C_0 \sqrt{N\log N}$), we would have that, for $k = N/\hat\mu_1+ \theta t_N /\hat \mu_2 $
\begin{align}
\hP&\Big( \htx_k= N , \hty_k= \gamma_c N +t_N, \text{``no big jump''}\Big)  =\frac{ 1+o(1)}{k} g \left( \frac{N-\hat \mu_1k}{\sqrt{k}}, \frac{\gamma_c N +t_N -\hat \mu_2 k}{\sqrt{k}} \right) \notag\\
& = \frac{(1+o(1))\hat \mu_1}{2\pi N \sqrt{(1-\rho^2) \sigma_1^2 \sigma_2^2}} \exp \bigg( -  \frac{\hat \mu_1 t_N^2}{2(1-\rho^2) N}  \Big\{ \frac{\theta^2}{\gamma_c^2 \sigma_1^2} - 2\rho \frac{\theta (1-\theta)}{\gamma_c \sigma_1\sigma_2} + \frac{(1-\theta)^2}{\sigma_2^2} \Big\} \bigg) \, ,
\label{localmoderate}
\end{align}
where $g(\cdot,\cdot)$ is the bivariate normal density of the limit $\frac{1}{\sqrt{k}} \big( \hat \tau_k - (\hat \mu_1, \hat \mu_2) k \big)$, which is centered with normalized covariance $\rho = (\sigma_1 \sigma_2)^{-1} {\rm Cov}(\htx_1,\hty_1) \in(-1,1)$ --- and $\sigma_1^2, \sigma_2^2$ are the respective variances of $\htx_1,\hty_1$. For the second equality, we used that $N-\hat \mu_1 k = \gamma_c^{-1} \theta t_N$, $\gamma_c N +t_N-\hat \mu_2 k = (1-\theta) t_N$ and $k = (1+o(1)) N/\hat \mu_1$.

Hence, in the sum over $k$ in \eqref{eq:conj1}, the main contribution should be for $k = \frac{N}{\hat\mu_1} + \frac{\theta_0}{\hat\mu_2} t_N   + O(\sqrt{N})$, with $\theta_0$ minimizing $Q(\theta):=\frac{\theta^2}{\gamma_c^2 \sigma_1^2} - 2\rho \frac{\theta (1-\theta)}{\gamma_c \sigma_1\sigma_2} + \frac{(1-\theta)^2}{\sigma_2^2}$ --- after some calculation we find that $\min Q(\theta) = (1-\rho^2) (\gamma_c^2 \sigma_1^2 + 2\rho \gamma_c \sigma_1 \sigma_2 +\sigma_2^2) ^{-1}$.
We end up with the following conjecture  in the case $\ga>1$, when $t_N/\sqrt{N} \to+\infty$
\begin{align}
\hP\Big( (N,\gamma_cN+t_N) \in\hat\tau\Big) =& (1+o(1)) \frac{N}{\hat \mu_1^2} \hP\Big( \hty_1=t_N\Big) + (1+o(1)) \frac{c_1}{\sqrt{N}} \exp\Big( -   \mathbf{c} \, \frac{t_N^2 }{N} \Big) \, ,
\label{eq:conj}
\end{align}
with $\mathbf{c}:= \frac{\hat \mu_1}{2}  (\gamma_c^2 \sigma_1^2 + 2\rho \gamma_c \sigma_1 \sigma_2 +\sigma_2^2) ^{-1}$, and the constant $c_1$ could in principle be made explicit.

\medskip
Plugging $t_N= a \sqrt{N\log N}$ in \eqref{eq:conj}, we find that the first term  is regularly varying with index $- \alpha/2 $ and that the second term has index $-1/2-  \mathbf{c} a^2 $.
Hence, depending on $a$ we can identify the dominant term in \eqref{eq:conj}:
\begin{equation}
\notag
 1^{\rm st} \text{ term is dominant }   \text{ if } a>  \sqrt{\frac{\ga-1}{2\mathbf{c}}} \, \quad \text{ and }  \quad
 2^{\rm nd} \text{ term is dominant }   \text{ if } a<  \sqrt{\frac{\ga-1}{2 \mathbf{c}}} \, .
\end{equation}
We therefore interpret this as  $a_c=  \sqrt{(\ga-1)/2\mathbf{c}} $, where $a_c$ is the critical value mentioned at the beginning of Appendix \ref{app2}, separating a big-jump domain (when $t_N >a \sqrt{N\log N}$ with $a>a_c$) from a moderate Cram\'er regime (when $t_N < a\sqrt{N\log N}$ with $a<a_c$).

\begin{figure}[htbp]
\centering
\includegraphics[width=14 cm]{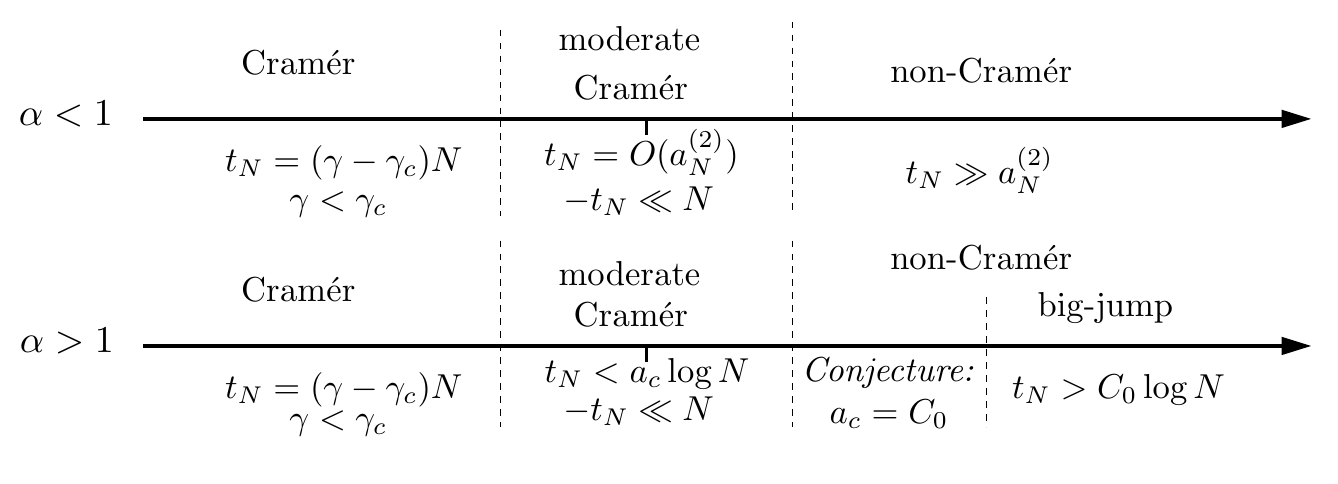}
\vskip-.4cm
\caption{\label{fig:flat} A schematic sum-up of the correspondence of the the values of $t_N%=M- \gamma_cN
= (\gamma-\gamma_c)N$
with the different regimes. We treat the big-jump domain and it  is the one to the right of the right-most dashed line. 
We believe that to the right of the moderate Cram\'er regime there is the big-jump domain -- put otherwise, the non-Cram\'er regime coincides
with the big-jump domain -- but this is proven only for $\ga<1$.}
\end{figure}

\smallskip
Notice that, when $t_N /\sqrt{N} \to -\infty$, one could develop an identical argument (except that the big-jump term disappears), provided that a local moderate deviation theorem  as \eqref{localmoderate} holds -- i.e. provided that $|t_N|/\sqrt{N}$ is not too large, how large depend mostly on the tail exponent $1+\ga>2$ of $\hat\tau^{(2)}$.
In the end, the sharp asymptotics of $\hP\big( (N,\gamma_cN+t_N) \in\hat\tau\big)$  should also be given by the second term in \eqref{eq:conj} -- as already seen in the case $t_N/\sqrt{N} \to t\in\bbR$ in Appendix \ref{app1}.

%In the case $\ga>1$, we have a full conjecture for the behavior of $\hP\Big( (N,\gamma_cN+t_N) \in\hat\tau\Big)$ on the critical Cram\'er region, that is when $t_N/N \to 0$.

\end{appendix}

%%%%%%%%%%%%%%%%%%%%%%%%%%%%%%%%%%%%%%%%%%%%%%%%%%%%%%%%%%%%%%%%%%%%%%%%%%%%%%
%%%%%%%%%%%%%%%%%%%%%%%%%%%%%%%%%%%%%%%%%%%%%%%%%%%%%%%%%%%%%%%%%%%%%%%%%%%%%%
%%%%%%%%%%%%%%%%%%%%%%%%%%%%%%%%%%%%%%%%%%%%%%%%%%%%%%%%%%%%%%%%%%%%%%%%%%%%%%
%%%%%%%%%%%%%%%%%%%%%%%%%%%%%%%%%%%%%%%%%%%%%%%%%%%%%%%%%%%%%%%%%%%%%%%%%%%%%%
%%%%%%%%%%%%%%%%%%%%%%%%%%%%%%%%%%%%%%%%%%%%%%%%%%%%%%%%%%%%%%%%%%%%%%%%%%%%%%

%These estimates lead to understanding the transitions and the different phases that one observe in the localized regime. We have seen that if $\gamma$ is in the Cram\'er regime, the free ends are microscopic, i.e. $O(1)$, and the limit process is just a recurrent renewal. By Proposition \ref{th:zfo}, if $\gamma$ in the interior of the complementary of the Cram\'er region instead a big loop appear and can be anywhere along the chain in the constrained case. For the free case, the excess of bases may be in the bulk ($\overline{\ga } > 1 + \ga$) or in the free ends ($\overline{\ga } < 1 + \ga$). 

\end{document}